\documentclass[12pt]{amsart}
\usepackage{latexsym,amsmath,amsfonts,amscd,amssymb}
\newtheorem{theorem}{\bf Theorem}[section]
\newtheorem{lemma}[theorem]{\bf Lemma}
\newtheorem{remark}[theorem]{\bf Remark}

\newtheorem{proposition}[theorem]{\bf Proposition}
\newtheorem{corollary}[theorem]{\bf Corollary}
\newtheorem{definition}[theorem]{\bf Definition}

\newtheorem{assumption}[theorem]{\bf Assumptions}
\newcommand{\be}{\begin{equation}}

\newcommand{\noi}{\noindent}

\newfont{\bfc}{cmbsy10 scaled 1200}  
\newfont{\dr}{msbm10 scaled \magstep1}  
\newfont{\sdr}{msbm8}  
\newfont{\gl}{eufm10 scaled \magstep1}  

\DeclareFontFamily{OT1}{rsfs}{}
\DeclareFontShape{OT1}{rsfs}{n}{it}{<->rsfs10}{}
\DeclareMathAlphabet{\curly}{OT1}{rsfs}{n}{it}

 \newcommand{\pf}{{\em Proof}. }

 \newcommand{\CC}{{\mathbb C}}

 \newcommand{\HH}{{\mathbb H}}
 
 \newcommand{\PP}{{\mathbb P}}
 
 \newcommand{\RR}{{\mathbb R}}

 \newcommand{\ZZ}{{\mathbb Z}}

 \newcommand{\cC}{{\mathcal C}}
 
 \newcommand{\cE}{{\mathcal E}}
 \newcommand{\cF}{{\mathcal F}}

 \newcommand{\cI}{{\mathcal I}}
 
 \newcommand{\cL}{{\mathcal L}}

 \newcommand{\cO}{{\mathcal O}}
 \newcommand{\cP}{{\mathcal P}}

 \newcommand{\cS}{{\mathcal S}}

 \newcommand{\cV}{{\mathcal V}}

\newcommand{\qu}{/\kern-.7ex/}
\newcommand{\exh}{\to\kern-1.8ex\to}

\newcommand{\End}{\operatorname{End}}
\newcommand{\ra}{\rightarrow}
\newcommand{\lra}{\longrightarrow}
\newcommand{\codim}{\operatorname{codim}}
\newcommand{\Pic}{\operatorname{Pic}}

\newcommand{\Hom}{\operatorname{Hom}}
\newcommand{\Ext}{\operatorname{Ext}}
\newcommand{\cExt}{\mathcal{E}xt\,}
\newcommand{\cHom}{\mathcal{H}om\,}
\newcommand{\cTor}{\mathcal{T}or\,}
\newcommand{\ndk}{(n,d,k)}

\newcommand{\Aut}{\operatorname{Aut}}
\newcommand{\Gr}{\operatorname{Gr}}
\newcommand{\rk}{\operatorname{rk}}
\newcommand{\id}{\mathrm{Id}}
\newcommand{\GCD}{\operatorname{GCD}}

\begin{document}
\title[On the geometry of moduli spaces of coherent systems]
{On the geometry of moduli spaces of coherent systems on algebraic curves}
\thanks{All authors are members of the research group VBAC
(Vector Bundles on Algebraic Curves), which was supported by
EAGER (EC FP5 Contract no.~HPRN-CT-2000-00099) and by EDGE (EC FP5
Contract no. HPRN-CT-2000-00101). Support was also received from
two grants from the European Scientific Exchange Programme of the
Royal Society of London and the Consejo
Superior de Investigaciones Cient\'{\i}ficas (15455 and 15646).
The first author was partially supported by the National Science
Foundation under grant DMS-0072073.}
\date{26 July 2006}
\author{S.~B.~Bradlow}
  \address{Department of Mathematics\\University of Illinois\\
  Urbana\\IL 61801\\USA}
  \email{bradlow@math.uiuc.edu}
\author{O.~Garc\'{\i}a-Prada}
   \address{Departamento de Matem{\'a}ticas \\
   Instituto de Matem\'aticas y F\'{\i}sica Fundamental\\
   Consejo Superior de Investigaciones Cient\'{\i}ficas \\ Serrano, 113 bis\\
   28006 Madrid \\ Spain}
   \email{oscar.garcia-prada@uam.es}
\author{V.~Mercat}
   \address{5 rue Delouvain, 75019 Paris, France}
   \email{mercat@math.jussieu.fr}
\author{V.~Mu\~noz}
   \address{Departamento de Matem{\'a}ticas \\
   Instituto de Matem\'aticas y F\'{\i}sica Fundamental\\
   Consejo Superior de Investigaciones Cient\'{\i}ficas \\ Serrano, 113 bis\\
   28006 Madrid \\ Spain}
  \email{vicente.munoz@imaff.cfmac.csic.es}
\author{P.~E.~Newstead}
   \address{Department of Mathematical Sciences \\
   University of Liverpool \\ Peach Street \\
   Liverpool L69 7ZL \\ UK}
   \email{newstead@liverpool.ac.uk}
\begin{abstract} Let $C$ be an algebraic curve of genus $g\ge2$.
A coherent system on $C$ consists of a pair $(E,V)$, where
$E$ is an algebraic vector bundle over $C$ of rank $n$ and
degree $d$ and $V$ is a subspace of dimension $k$ of the
space of sections of $E$. The stability of the coherent
system depends on a parameter $\alpha$.
We study the geometry of the moduli space of coherent
systems for different values of $\alpha$ when $k\leq n$ and
the variation of the moduli spaces when we vary $\alpha$.
As a consequence, for sufficiently large $\alpha$,
we compute the Picard groups and the first and second
homotopy groups of the moduli spaces of coherent systems
in almost all cases, describe the moduli space for the case $k=n-1$
explicitly, and give the Poincar\'e polynomials
for the case $k=n-2$. In an appendix, we describe the geometry of the ``flips'' 
which take place at critical values of $\alpha$ in the simplest case, and
include a proof of the existence of universal families of coherent systems
when $\GCD(n,d,k)=1$.

\end{abstract}
\maketitle

\begin{center}{\em Dedicated to the memory of Joseph Le Potier}\end{center}

\section{Introduction}\label{section:intro}
Let $C$ be a smooth projective algebraic curve of genus $g\geq 2$. A coherent
system on $C$
of type $(n,d,k)$ is a pair $(E,V)$, where $E$ is a vector bundle on $C$ of rank
$n$ and
degree $d$ and $V$ is a subspace of dimension $k$ of the space of sections
$H^0(E)$.
Introduced in \cite{KN}, \cite{RV} and \cite{LeP}, there is a notion of
stability for
coherent systems which permits the construction of moduli spaces. This notion
depends on a
real parameter, and thus leads to a family of moduli spaces. As described in
\cite{BG},
there is a useful relation between these moduli spaces and the Brill-Noether
loci in the
moduli spaces of semistable bundles of rank $n$ and degree $d$.

In \cite{BGMN} we began a systematic study of the coherent systems moduli
spaces, partly
with a view to applications in higher rank Brill-Noether theory.  In this paper,
we continue
our explorations and obtain substantial new information about the geometry and
topology of
the moduli spaces in the case $k\le n$. In particular we obtain precise
conditions for the moduli space to be non-empty and show that each non-empty
moduli space has a distinguished irreducible component which has the expected
dimension. For sufficiently large values of the parameter, we show that this
is the only component. In addition, we obtain some precise information on
the way in which the moduli space changes as the parameter varies, and deduce
from this some
more refined results on Picard groups and
Picard varieties and on homotopy and cohomology groups.
When $k=n-1$, we obtain a complete geometrical description of the moduli space
for sufficiently large values of the parameter; this allows us to give a
precise description of certain Brill-Noether loci. This is the only
application to Brill-Noether theory which we include here; our results will
certainly make further contributions to this theory, but this involves
additional technicalities and will be addressed in future papers.

We refer the reader to \cite{BGMN} (and the further references cited therein)
for the basic
properties of coherent systems on algebraic curves. For convenience, and in
order to set
notation, we give a short synopsis before proceeding to a more detailed
description of our results.

\begin{definition} \label{def:stable}
{\em Fix $\alpha \in \RR$. Let $( E,V)$\ be a coherent system of type $\ndk$. The
$\alpha$-{\em
slope} $\mu_{\alpha}(E,V)$ is defined by
 $$
 \mu_{\alpha}( E,V)=\frac{d}{n}+\alpha\frac{k}{n}\ .
 $$
We say $( E,V)$\ is $\alpha$-{\em stable} if
 $$
 \mu_{\alpha}(E',V')<\mu_{\alpha}(E,V)
 $$
for all proper subsystems $(E',V')$ \textup{(}i.e. for every non-zero subbundle
$E'$ of $E$
and every subspace $V'\subseteq V\cap H^0(E')$ with $(E',V')\ne(E,V)$\textup{)}.
We define
$\alpha$-{\em semistability} by replacing the above strict inequality with a
weak
inequality. A coherent system is called $\alpha$-{\em polystable} if it is the
direct sum of
$\alpha$-stable coherent systems of the same $\alpha$-slope. We denote the {\em
moduli
space} of $\alpha$-stable coherent systems of type $\ndk$ by $G(\alpha;n,d,k)$.}

\end{definition}

\begin{definition}\label{def:cv}~ {\em

\begin{itemize}

\item We say that  $\alpha>0$\  is a  {\em critical value} (or, in the
terminology of \cite[Definition 2.4]{BGMN}, {\em actual critical value}) if
there exists a proper subsystem
$(E',V')$\ such that $\frac{k'}{n'}\ne\frac{k}{n}$\ but $\mu_{\alpha}(
E',V')=\mu_{\alpha}(
E,V)$. We also regard $0$ as a critical value.

\item We say that $\alpha$ is {\em generic} if it is not a critical value.
If $\GCD(n,d,k)=1$ and
$\alpha$ is generic, then $\alpha$-semistability is equivalent to $\alpha$-
stability.

\item If we label the critical values of $\alpha$\ by $\alpha_i$, starting with
$\alpha_0=0$,
we get a partition of the $\alpha$-range into a set of intervals
$(\alpha_{i},\alpha_{i+1})$. Within the interval $(\alpha_{i},\alpha_{i+1})$ the
property of
$\alpha$-stability is independent of $\alpha$, that is if $\alpha,\alpha'\in
(\alpha_{i},\alpha_{i+1})$, $G(\alpha;n,d,k)= G(\alpha';n,d,k)$. We shall denote
this moduli
space by $G_i=G_i(n,d,k)$.
\end{itemize}}
\end{definition}

The construction of moduli spaces thus yields one moduli space $G_i$ for the
interval
$(\alpha_i,\alpha_{i+1})$. If $\GCD(n,d,k)\ne 1$, one can define similarly the
moduli spaces
$\tilde{G}_i$ of semistable coherent systems. The GIT construction of these
moduli spaces
has been given in \cite{LeP} and \cite{KN}. A previous construction for $G_0$
had been given
in \cite{RV} and also in \cite{Be} for large values of $d$.

In the next two propositions, we suppose that $G(\alpha;n,d,k)\ne\emptyset$
for at least one value of $\alpha$.

\begin{proposition}{\bf(\cite[Proposition 4.2]{BGMN})}\label{prop:last}
Let $0<k<n$ and let $\alpha_L$ be  the biggest critical value smaller than
$\frac{d}{n-k}$.
The $\alpha$-range is then divided into a finite set of intervals
determined by critical values
 $$
 0=\alpha_0<\alpha_1<\dots <\alpha_L<\frac{d}{n-k}\ .
 $$
If
$\alpha> \frac{d}{n-k}$, the moduli spaces are empty.
\end{proposition}

We shall also be concerned with the case $k=n$. For this we need the following
result from \cite{BGMN}.

\begin{proposition}{\bf(\cite[Proposition 4.6]{BGMN})}\label{prop:flips}
Let $k\geq n$. Then there is a critical value, denoted by $\alpha_L$, after
which the moduli
spaces stabilise, i.e. $G(\alpha;n,d,k)=G_L$ if $\alpha>\alpha_L$. The $\alpha$-range is
thus divided into a finite set of intervals bounded by critical values
 $$
 0=\alpha_0<\alpha_1<\dots<\alpha_L<\infty
 $$
such that,
for any two different values of $\alpha$ in the range $(\alpha_L,\infty)$, the
moduli spaces
coincide.
\end{proposition}

The difference between adjacent moduli spaces in the family $G_0,G_1,\dots,G_L$
is accounted
for by the subschemes $G_i^+\subseteq G_i$, $G_i^-\subseteq G_{i-1}$  where $G_i^+$ consists of all
$(E,V)$ in
$G_i$ which are not $\alpha$-stable if $\alpha<\alpha_i$ and $G_i^-\subseteq
G_{i-1}$
contains all $(E,V)$ in $G_{i-1}$ which are not $\alpha$-stable if
$\alpha>\alpha_i$. It
follows that $G_i-G_i^+$ and $G_{i-1}-G_i^-$ are isomorphic and that $G_i$ is
transformed
into $G_{i-1}$ by the removal of $G_i^+$ and the insertion of $G_i^-$.

\begin{definition}
{\em We refer to such a procedure, i.e. the transformation of $G_i$ into
$G_{i-1}$ by the removal of $G_i^+$ and the insertion of $G_i^-$, or the
inverse transformation, as a {\em flip}. We say
that a flip is {\em good} if the flip loci (i.e. the subschemes $G_i^{\pm}$)
have strictly
positive codimension in every component of the moduli spaces $G_i$, $G_{i-1}$; under these conditions
the
moduli spaces $G_i$, $G_{i-1}$ are birationally equivalent.}
\end{definition}

Our main results fall into three categories: non-emptiness, smoothness
and irreducibility of the
moduli spaces; homotopy and Picard groups and Picard varieties;
Poincar\'e polynomials.

In \cite{BG} and \cite{BGMN} we obtained non-emptiness and irreducibility
results for the
moduli space $G_L(n,d,k)$. In this paper we extend these results
to other moduli spaces $G(\alpha;n,d,k)$. Our results are summarized in the
following theorems.

\noindent{\bf Theorem A [Lemmas \ref{lemma:inj2} and \ref{lemma:useful}, and
 Theorems \ref{th:inj} and \ref{th:alphai}]}
{\it Suppose that $0<k\le n$ and $n\ge2$. Then the moduli space
$G(\alpha;n,d,k)$ is non-empty if and only if

$$\alpha>0\ ,\ (n-k)\alpha<d\ ,\ k\le n+\frac{1}{g}(d-n)\ ,\ \mbox{and}\
(n,d,k)\ne(n,n,n) \ .$$

Whenever it is non-empty, $G(\alpha;n,d,k)$ contains a Zariski-open subset,
$U(\alpha)$, which is smooth, irreducible of dimension
$$
\beta(n,d,k)=n^2(g-1)+1-k(k-d+n(g-1)),
$$ and
such that its closure $\overline{U(\alpha)}\subset G(\alpha;n,d,k)$ is
birationally equivalent to $G_L$.

There is a critical value $\alpha_I\in [0,\frac{d}{n-k})$ such that for all
$\alpha>\alpha_I$ the moduli space $G(\alpha;n,d,k)=U(\alpha)$. Thus for
$\alpha>\alpha_I$, whenever it is non-empty, $G(\alpha;n,d,k)$ is smooth,
irreducible of dimension $\beta(n,d,k)$ and birationally equivalent to
$G_L(n,d,k)$. The critical value $\alpha_I$ satisfies the bound

$$\alpha_I\le \max\left\{\frac{(k-1)(d-n)-n\epsilon}{k(n-k+1)},0\right\}$$

\noindent where $\epsilon=\min\{k-1,g\}$. }

One consequence of the above theorem is that, if $k\le n$,
then there are precisely two possibilities for the
moduli spaces $G(\alpha;n,d,k)$: for fixed $(n,d,k)$, either
$G(\alpha;n,d,k)$ is non-empty for all allowable $\alpha$, or it is empty
for all $\alpha$. Moreover, the non-empty
moduli spaces always contain a distinguished component of the expected
dimension (i.e. $\beta(n,d,k)$) -- we will identify this component more
precisely
in section \ref{section:galpha}.

For the sake of completeness, we recall from  \cite{BGMN}  that if $k< n$,
then, whenever it is non-empty, the moduli space $G_L(n,d,k)$ is birationally
equivalent
to a fibration over $M(n-k,d)$, the moduli space of stable bundles of rank
$n-k$ and degree $d$, with fibre the Grassmannian $\Gr(k,d+(n-k)(g-1))$. If
$\GCD(n-k,d)=1$, then the birational equivalence is an isomorphism.

Our next result concerns the case $k=n-1$.

\noindent{\bf Theorem B [Corollary \ref{cor:k=n-1} and Theorem \ref{th:gt}]}
{\it Suppose that
$n\ge2$ and $d>0$. For all $\alpha$ such that
$$\max\{d-n,0\}<\alpha<d$$

\noindent the moduli space $G(\alpha;n,d,n-1)$ is non-empty if and only if
$$d\ge\max\{1,n-g\}\ .$$
\noindent Moreover, whenever it is non-empty, $G(\alpha;n,d,n-1)=G_L(n,d,n-1)$
and is a
fibration over the Jacobian $J^d$, with fibre the Grassmannian
$\Gr(n-1,d+g-1)$.} (This fibration will be identified more precisely in
section \ref{section:n-1}.)

As a consequence of Theorem B, we identify the Brill-Noether locus
$B(n,d,n-1)$, consisting of stable bundles $E$ of rank $n$ and degree $d$ with
$h^0(E)\ge n-1$, with a well known classical variety (Theorem \ref{th:bnn-1}).

Our main results on the Picard groups and varieties and homotopy groups for $k<n$
are
summarized in the following theorem.

\noindent{\bf Theorem C [Theorem \ref{th:gt}, Theorem \ref{picard},
Corollary \ref{cor:picard}
and Theorems \ref{homotopy:k=n-1}, \ref{homotopy}, \ref{picard-general}
and \ref{picvar}]}
{\it Let $0<k<n$ and $d>0$. Suppose further that $k<n+\frac1g(d-n)$ and
$\max\left\{\frac{d-n}{n-k},0\right\}<\alpha<\frac{d}{n-k}$. Then, except possibly when $g=2$,
$k=n-2$ and $d$ is even,

\begin{itemize}

\item[(a)] $ \Pic(G(\alpha;n,d,k))\cong \Pic(M(n-k,d))\times
\ZZ$;

\item[(b)] $\Pic^0(G(\alpha;n,d,k))$ is isomorphic to the Jacobian $J(C)$;
\item[(c)] $\pi_1(G(\alpha;n,d,k))\cong
\pi_1(M(n-k,d))\cong H_1(C,\ZZ)$;

\item[(d)] if also $k\ne n-1$, there exists an exact sequence 
$$0\lra\ZZ\lra\pi_2(G(\alpha;n,d,k)) \lra\ZZ \times\ZZ_p\lra0,$$ where $p=\GCD(n-k,d)$.
\end{itemize}
If $k=n+\frac1g(d-n)$, {\em (b)} and {\em (c)} remain true; in {\em(a)},
the factor $\ZZ$ must be deleted; in {\em (d)}, $\pi_2(G(\alpha;n,d,k)) \cong\ZZ \times\ZZ_p$.

In the case $k=n-1$, we have
\begin{itemize}
\item [(e)] if $\max\{d-n,0\}<\alpha<d$ and  $d\geq \max\{1,n-g\}$, then
$$\pi_i (G(\alpha;n,d,n-1)) \cong
\pi_i(\Gr(n-1,d+g-1))\ \mbox{for}\ i\geq 0, i\ne1.$$
\end{itemize}}

Finally we have some results on Poincar\'e polynomials in the case $k=n-2$.
For any space $X$, we write $P(X)$ for the Poincar\'e polynomial of $X$
with coefficients in any fixed field.
In order to state our results concisely, we write  $G(\alpha)$, $G_L$ for
$G(\alpha;n,d,n-2)$, $G_L(n,d,n-2)$. We write also, for any $r$, $e$,
$$
P_{r,e}=P(G_L(r,e,r-1)).
$$

\noindent{\bf Theorem D [Theorem \ref{poincare2}, Corollary \ref{cor:n=3},
Corollary \ref{cor:n=4}]}
{\it
\begin{itemize}
\item[(a)]For any
non-critical value $\alpha'$ in the interval
$(\max\{\frac{d-n}2,0\},\frac{d}2)$,
$$
 P(G(\alpha'))(t)-P(G_L)(t)=\sum\frac{t^{2C_{21}(n_1,d_1)}-
 t^{2C_{12}(n_1,d_1)}}{1-t^2}P_{n_1,d_1}(t)P_{n-n_1,d-d_1}(t),
$$
where the summation is over all solutions of {\rm(\ref{eq:n1})},
{\rm(\ref{eq:alphan-2})}, {\rm(\ref{eq:d1})} for which $\alpha>\alpha'$.
{\rm(For the definitions of the polynomials $C_{12}$ and $C_{21}$ and the
equations  (\ref{eq:n1}),
(\ref{eq:alphan-2}), (\ref{eq:d1}), see section \ref{section:n-2}.)}

\item[(b)]Suppose $n=3$ and $d$ is odd. Then, for
$\max\{\frac{d-3}2,0\}<\alpha'<\frac{d}2$,
\begin{eqnarray*}
P(G(\alpha'))(t)&=&P(G_L)(t)=P(M(2,d))(t)\frac{1-t^{2(d+2g-2)}}{1-t^2}\\
&=&\frac{(1+t)^{2g}((1+t^3)^{2g}-t^{2g}(1+t)^{2g})(1-t^{2(d+2g-2)})}{(1-
t^2)^2(1-t^4)}.
\end{eqnarray*}

\item[(c)]Suppose $n=4$ and $d$ is odd. Then
\begin{itemize}
\item[\rm(i)] if $\max\{\frac{d-2}2,0\}<\alpha'<\frac{d}2$,
\begin{eqnarray*}
P(G(\alpha'))(t)&=&P(G_L)(t)\\&=&
\frac{(1+t)^{2g}((1+t^3)^{2g}-t^{2g}(1+t)^{2g})(1-t^{2(d+2g-3)})
(1-t^{2(d+2g-2)})}{(1-t^2)^2(1-t^4)^2};
\end{eqnarray*}
\item[\rm(ii)] if $\max\{\frac{d-4}2,0\}<\alpha'<\frac{d-2}2$,
\begin{eqnarray*}
P(G(\alpha'))(t)
&=&\frac{(1+t)^{2g}((1+t^3)^{2g}-t^{2g}(1+t)^{2g})
(1-t^{2(d+2g-3)})(1-t^{2(d+2g-2)})}{(1-t^2)^2(1-t^4)^2}\\
&&+\frac{(t^{2g}-t^{6g+2d-10})(1-t^{d-3+2g})
(1-t^{d-1+2g})(1+t)^{4g}}{(1-t^2)^2(1-t^4)}.
\end{eqnarray*}
\end{itemize}
\end{itemize}}

 We now summarize the contents of the paper. In section \ref{section:klen}, we
show that,
for $k\le n$ and $\alpha$ sufficiently large, any $\alpha$-stable coherent
system $(E,V)$ is
injective, in other words $(E,V)$ gives rise to an injective morphism
$V\otimes {\cO}\hookrightarrow E$. In section \ref{section:galpha} we
obtain precise conditions for non-emptiness of the moduli spaces of
injective coherent systems and deduce our first irreducibility
results. In section \ref{section:diag}
we assume $k<n$ and prove an important structural result (the Diagram Lemma)
for torsion-free coherent systems (that is, injective coherent systems with
$E/(V\otimes\cO)$ torsion-free).
In section \ref{section:n-1} we give a first application of this lemma
to the moduli spaces of coherent systems with $k=n-1$ and to the
Brill-Noether locus $B(n,d,n-1)$. Following this, in section
\ref{section:flip} we compute codimensions for the flips
which occur at critical values
$\alpha>\frac{d-n}{n-k}$ when $k<n$. In section \ref{section:picard} we obtain
results on
Picard groups and varieties and homotopy groups. In section \ref{section:n-2} we
apply the calculations of
section \ref{section:flip} to the case $k=n-2$ and investigate the geometry of
the flips in
a way similar to the work of Thaddeus; this allows a computation of Poincar\'e
polynomials. The necessary extensions to Thaddeus' results are contained in an 
appendix. We give also in the appendix a proof of the existence of universal 
families of coherent systems when $\GCD(n,d,k)=1$, as we were unable to locate
a proof in the literature.

We suppose throughout that $C$ is a smooth irreducible projective algebraic
curve of genus
$g\ge2$ defined over the complex numbers. The special cases $g=0$ and $g=1$ are
being investigated elsewhere
\cite{LN1,LN2,LN3}.

\section{Coherent Systems with $k\leq n$}\label{section:klen}

Let $(E,V)$ be a coherent system of some fixed type $(n,d,k)$ on
$C$ with $k\ge1$. For most of this section we assume that $k\leq
n$. For convenience we introduce the following definition.

\begin{definition}\label{def:inj}
{\em A coherent system $(E,V)$ is {\em injective} if the evaluation
morphism $V\otimes\cO\ra E$ is injective as a morphism of sheaves.
Moreover $(E,V)$ is {\em torsion-free} if it is injective and the
quotient sheaf $E/(V\otimes\cO)$ is torsion-free.}
\end{definition}

Since we are working over a smooth curve, we have for any
torsion-free coherent system $(E,V)$ an exact sequence
\be\label{eq:tf} 0\lra V\otimes\cO\lra E\lra F\lra0,
\end{equation}
where $F$ is a vector bundle on $C$.

\begin{lemma}\label{lemma:inj}
Suppose $(E,V)$ is injective and $\alpha$-stable. Then
$G(\alpha;n,d,k)$ is smooth of dimension
 \be\label{eq:beta}
 \beta(n,d,k)=n^2(g-1)+1-k(k-d+n(g-1))
 \end{equation}
at $(E,V)$.
\end{lemma}

\pf This is \cite[Proposition 3.12]{BGMN}.\qed

We recall from \cite[Lemma 3.6]{BG} and \cite[Proposition
4.4]{BGMN} that, for $k\le n$, every $\alpha$-semistable coherent
system is injective provided $\alpha$ is large enough. Moreover,
if $k<n$, then, by \cite[Lemma 3.7]{BG} every $\alpha$-semistable
coherent system is torsion-free, again for $\alpha$ large enough.

\begin{definition}\label{def:alpha}
{\em We define $\alpha_I\ge0$ to be the smallest critical value of
$\alpha$ such that every $\alpha$-semistable coherent system is
injective  for $\alpha> \alpha_I$. If $k<n$, we define
$\alpha_T\ge0$ to be the smallest critical value of $\alpha$ such
that every $\alpha$-semistable coherent system is torsion-free for
$\alpha> \alpha_T$.}
\end{definition}

Of course, $\alpha_I$ and $\alpha_T$ depend on the type $(n,d,k)$
and we have $\alpha_T\ge\alpha_I$. Note that, if $k=1$,
$\alpha_I=0$. We have

\begin{lemma} \label{lemma:inj2}Let $k\le n$. Then
\begin{itemize}
\item[(i)]\begin{itemize} \item[(a)] if $2\le k\leq
g+1$, $\alpha_I\leq
\max\left\{\frac{(k-1)(d-2n)}{k(n-k+1)},0\right\}$;
\item[(b)] if $k>g+1$,
$\alpha_I\leq\max\left\{\frac{(k-1)(d-n)-ng}{k(n-k+1)},0\right\}$;
\end{itemize}
\item[(ii)] if $0<k<n$, $\alpha_T=
\max\left\{\frac{d-n}{n-k},0\right\}$.
\end{itemize}
\end{lemma}
\pf (i) Suppose the evaluation morphism is not injective. Consider
the subsheaf $Im$  generated by $V$ in $E$. It is a standard fact
that $Im\simeq F\oplus {\cO}^s$, where $s=h^0(Im^*)$ and $F$ is a
vector bundle. Clearly $h^0(F^*)=0$ and $F$ is generated by the
global sections which lie in the image $W'$  of $V$ in $H^0(F)$. We
write $l=\rk F$ and $e=\deg F$.

It is well known that a general subspace $W\subset W'$ of
dimension $l+1$ generates $F$. We have therefore an exact sequence
$$
0\longrightarrow (\det F)^*\longrightarrow W\otimes {\cO}\longrightarrow F
\longrightarrow 0\ ,
$$
and by dualising we obtain that $h^0(\det F)\geq l+1$. It follows
easily from Clifford's Theorem and the Riemann-Roch Theorem that
\begin{itemize}
\item[($*$)] if $l\leq g$, then $e\geq 2l$; \item[($**$)] if
$l\geq g$, then $e\geq l+g$.
\end{itemize}
Note that by definition we have $l\leq k-1$.

Now the $\alpha$-semistability criterion implies that
$\mu_\alpha(F,W)\le\mu_\alpha(E,V)$ and hence
$$\frac{e+\alpha (l+1)}{ l}\le\frac{d+\alpha k}{ n},$$
or equivalently
$$\alpha(n(l+1)-kl)\le ld-ne.$$
Since $n\geq k$, we necessarily have $n(l+1)-kl>0$ and the above
inequality gives
$$\alpha\le\frac{ld-ne}{n(l+1)-kl}\ .$$
 Using the inequalities $(*)$, $(**)$ and the condition $l\leq k-1$, we obtain
(i).\\

(ii) Suppose $0<k<n$. To show that $\alpha_T\ge
\max\left\{\frac{d-n}{n-k},0\right\}$, we need to show that, if
$d>n$, there exists an $\alpha$-semistable coherent system with
$\alpha=\frac{d-n}{n-k}$ which is not torsion-free. For this we
can take
$$
\bigoplus_{i=1}^k(\cO(p_i),H^0(\cO(p_i)))\oplus(F,0),
$$
where $p_1,\ldots,p_k\in C$ and $F$ is a semistable bundle of rank $n-k$ and
degree $d-k$.

To prove the opposite inequality, let $(E,V)$ be an
$\alpha$-semistable coherent system which is not torsion-free.
Either $(E,V)$ is not injective or the quotient sheaf
$E/(V\otimes\cO)$ has non-zero torsion. In both cases, there
exists a non-zero section $s$ in $V$ that vanishes at some point
of $C$. The section $s$ generically generates a line subbundle $L$
with $\deg L>0$. We thus have a coherent subsystem $(L,V_L)$ with
$\deg L\ge 1$, $\rk L=1$, and $\dim V_L\ge 1$.  The
$\alpha$-semistability criterion now yields

$$
1+\alpha\le\frac{d}{n}+\alpha\frac{k}{n},
$$
i. e.
$$
\alpha\le\frac{d-n}{n-k}.
$$
\qed

\begin{corollary}\label{cor:inj} If $\alpha>0$, $2\le k\le n$ and
$$d
\le\min\left\{2n,n+\frac{ng}{k-1}\right\},
$$
then every $\alpha$-semistable coherent system $(E,V)$ is
injective.
\end{corollary}

\pf This follows immediately from Lemma \ref{lemma:inj2}.\qed
\bigskip

\begin{remark}\label{rmk1}
{\em (i) It is a simple exercise to show that, if $k<n$ and $d>n$,
the bounds on $\alpha_I$ of Lemma \ref{lemma:inj2}(i) are always
strictly smaller than the value of $\alpha_T$ given by Lemma
\ref{lemma:inj2}(ii). In particular, if $\alpha_T>0$, then
$\alpha_I<\alpha_T$.

(ii) Note that, for $k<n$ and $d>0$, $\alpha_T<\frac{d}{n-k}$,
which is the maximum value of $\alpha$ for which
$\alpha$-semistable coherent systems exist (see \cite[Lemma
4.1]{BGMN}). Moreover, it was proved in \cite[Lemmas 3.6 and
3.7]{BG} that
$$
\alpha_I\leq \frac{d(k-1)}{k(n-k+1)},\ \
\alpha_T\leq\max\left\{\frac{kd-n}{k(n-k)},0\right\}.
$$
The statements of Lemma \ref{lemma:inj2} are stronger.

(iii) Compare also \cite[Chapitre 3, Lemme A.2]{M1}, in which it
is proved that, if $d<\min\{2n,n+g\}$ and $E$ is a semistable
bundle, then $(E,V)$ is injective. When $E$ is stable, Corollary
\ref{cor:inj} gives the same result with a weaker restriction on
$d$.}
\end{remark}

\begin{remark}\label{rmk4}
{\em In proving that, when $d>n$, $\alpha_T\ge\frac{d-n}{n-k}$, we
have made use of a coherent system which is $\alpha$-semistable
for $\alpha=\frac{d-n}{n-k}$ only. In fact it is possible to find
a coherent system $(E,V)$ which is not torsion-free and is
$\alpha$-stable for $\alpha$ slightly less than $\frac{d-n}{n-k}$.
For this, one can take $(E,V)$ to be given by a non-trivial
extension
$$
0\lra(E_1,V_1)\lra(E,V)\lra(E_2,V_2)\lra0,
$$
where $(E_1,V_1)=(\cO(p),H^0(\cO(p)))$ and $(E_2,V_2)$ is
$\alpha$-stable for $\alpha=\frac{d-n}{n-k}$. For the existence of
$(E_2,V_2)$, see Theorem \ref{th:inj} below and, for the existence
of a non-trivial extension, see equations (\ref{C21}) and
(\ref{dim-ext}).}
\end{remark}

\begin{remark}\label{rmk9}
{\em For a general curve, the bound in Lemma
\ref{lemma:inj2}(i)(a) is not best possible when $k<g$; we can
improve the bound by using an estimate based on the Brill-Noether
number rather than Clifford's Theorem. Details are left to the
reader.}
\end{remark}
The following lemma is true without the restriction $k\le n$.

\begin{lemma}\label{lemma:dualcs}
Let $(E,V)$ be an $\alpha$-stable coherent system for some
$\alpha>0$ with $k>0$. Then
\begin{itemize}
\item $d>0$, except in the case $(n,d,k)=(1,0,1)$. \item $h^0(E^*)=0$.
\end{itemize}
If $(E,V)$ is $\alpha$-semistable, then
\begin{itemize}
\item $d\ge0$, \item $h^0(E^*)=0$ except when $d=\alpha(n-k)$.
\end{itemize}
\end{lemma}

\pf For the fact that $d>0$ ($d\ge0$) for an $\alpha$-stable
($\alpha$-semistable) coherent system, see \cite[Lemmas 4.1 and
4.3]{BGMN}. Now suppose $(E,V)$ is $\alpha$-semistable and
$h^0(E^*)\ne0$. Then there exists a non-zero homomorphism
$E\ra\cO$. If the induced map $V\otimes {\mathcal O}\rightarrow
{\cO}$ is not the zero map, then $(\cO,H^0(\cO))$ is a direct
factor of $(E,V)$. This contradicts $\alpha$-stability always and
$\alpha$-semistability unless
$$
\alpha=\frac{d}{n}+\alpha\frac{k}{n},
$$
i.~e. $d=\alpha(n-k)$. Otherwise let $E'$ be the kernel of $E\ra\cO$. Then $(E',V)$ is a
coherent subsystem of $(E,V)$ with $\deg E'\ge d$, and the
$\alpha$-semistability criterion gives
$$
\frac{d}{n-1}+\alpha\frac{k}{n-1}\le\frac{d}{n}+\alpha \frac{k}{n}.
$$
This contradicts the assumption that $\alpha>0$. \qed

Returning now to the case $k\le n$, we have

\begin{lemma}\label{lemma:extclass}
Suppose that $k\le n$ and $d>0$.  Let $(E,V)$\ be any injective
coherent system of type $(n,d,k)$, i.e. suppose that $(E,V)$ is
represented by an extension
\begin{equation}\label{sequence:inj}
0\longrightarrow\cO^{\oplus k}\longrightarrow E\longrightarrow
F\longrightarrow 0
\end{equation}
(where $ F$\ need not be locally free).  Let
$$
\vec{e}=(e_1,\dots,e_k)\in \Ext^1(F,\cO^{\oplus k}) =\Ext^1(
F,\cO)^{\oplus k}
$$
denote the extension class of {\em(\ref{sequence:inj})}. If
$h^0(E^*)=0$, then
\begin{itemize}
\item $e_1,\ldots,e_k$ are linearly independent as vectors in
$\Ext^1( F,\cO)$; \item $h^0(F^*)=0.$
\end{itemize}
Moreover
\begin{equation}\label{eq:kbound}
k\le n+\frac1g(d-n).
\end{equation}
\end{lemma}

\pf Suppose that $e_1,\ldots,e_k$ are linearly dependent. After
acting on (\ref{sequence:inj}) by an automorphism of $\cO^{\oplus
k}$, we can suppose some $e_i=0$. But then $\cO$ is a direct
factor of $E$ and $h^0(E^*)\ne0$. This proves the first statement.

The vanishing of $h^0(F^*)$ follows from the long exact sequence
of the $\Ext^i(\cdot,\cO)$ induced by (\ref{sequence:inj}). By
Riemann-Roch, this implies that
$$
\dim\Ext^1( F,\cO)=d+(n-k)(g-1)\ .
$$
The bound in (\ref{eq:kbound}) follows from this and the linear
independence of $e_1,\dots,e_k$. \qed

\section{The moduli space for injective coherent systems}\label{section:galpha}

In this section we show that the moduli space of injective
$\alpha$-stable coherent systems is always smooth and irreducible
of dimension $\beta(n,d,k)$ whenever it is non-empty. We also
determine exactly when this space is non-empty. As consequences,
we obtain the same properties for $G(\alpha;n,d,k)$ when
$\alpha>\alpha_I$ and results on the birational type of some
moduli spaces. These results can be seen as an extension of
\cite[Theorem 5.4]{BGMN}.

We begin with a proposition which restates some key results from
\cite{BG} and \cite{BGMN} and a useful lemma.

\begin{proposition}\label{prop:restate}
{\em(i)} Suppose $n\ge2$ and $0<k\le n$. Then
$G_L(n,d,k)\ne\emptyset$ if and only if \be\label{eq:cond2} d>0,\
k\le n+\frac1g(d-n)\ \mbox{and}\ (n,d,k)\ne(n,n,n),
\end{equation}
and it is then always irreducible and smooth of dimension
$\beta(n,d,k)$.

{\em(ii)} If $k=n$, every element of $G_L(n,d,k)$ can be
represented by an extension of the form \be\label{eq:tor1} 0\lra
V\otimes\cO\lra E\lra T\lra 0,
\end{equation}
where $T$ is a torsion sheaf.

{\em(iii)} If $0<k<n$, every element of $G_L(n,d,k)$ is
torsion-free and corresponds to an extension {\em(\ref{eq:tf})} with $F$
semistable
and $h^0(F^*)=0$; moreover $G_L(n,d,k)$ is birationally equivalent to
a fibration over the moduli space $M(n-k,d)$ with fibre the
Grassmannian $\Gr(k,d+(n-k)(g-1))$. More precisely, if $W$ denotes the
subvariety of $G_L(n,d,k)$ consisting of coherent systems for which the
bundle $F$ in {\em(\ref{eq:tf})} is strictly semistable, then
$G_L(n,d,k)\setminus W$ is isomorphic to a Grassmann fibration over
$M(n-k,d)$.

{\em(iv)} If in addition $\GCD(n-k,d)=1$, then $W=\emptyset$ and
$G_L(n,d,k)\to M(n-k,d)$ is the Grassmann fibration associated to some
vector bundle over $M(n-k,d)$.
\end{proposition}
\pf For the necessity of the condition $d>0$, see Lemma
\ref{lemma:dualcs}. The rest is essentially a restatement of \cite[Proposition
5.2 and
Theorem 5.4]{BGMN} (see also \cite[section 4]{BG}) and \cite[Theorem
5.6]{BGMN}. (For the last part of (iv), note that there exists a Poincar\'e
bundle $\cP$ over $X\times M(n-k,d)$. The sheaf $R^1p_{2*}{\cP}^*$ is
locally free since $h^1(F^*)$ is constant for $F\in M(n-k,d)$.
The vector bundle corresponding to this sheaf has the required properties.)
\qed

\begin{lemma}\label{lemma:useful} Suppose $0<k\le n$ and let
$$
U=\{(E,V)\in G_L(n,d,k)\,|\,E\ \mbox{is stable}\}.
$$
Then $U\subset G(\alpha;n,d,k)$ is a Zariski-open subset for all
allowable $\alpha$ (that is, for $0<\alpha<\frac{d}{n-k}$ if $k<n$
and for $\alpha>0$ if $k=n$).
\end{lemma}
\pf If $(E,V)\in U$, then $(E,V)$ is $\alpha$-stable for small
$\alpha$ and for large $\alpha$. Since the set of $\alpha$ for
which $(E,V)$ is $\alpha$-stable is an open interval \cite[Lemma
3.14]{BG}, it follows that $(E,V)$ is $\alpha$-stable for all
allowable $\alpha$. This proves that $U\subset G(\alpha;n,d,k)$.
The fact that $U$ is Zariski-open follows from the openness of the
stability condition. \qed

We come now to the main result of this section.
\begin{theorem}\label{th:inj} Suppose $0<k\le n$. For any $\alpha>0$, define
$U(\alpha)$ by
$$
U(\alpha)=\{(E,V)\in G(\alpha;n,d,k)\,|\,(E,V)\ \mbox{is
injective}\}.
$$
Then
\begin{itemize}
\item[(i)] $U(\alpha)$ is a Zariski-open subset of
$G(\alpha;n,d,k)$;
\item[(ii)] if $U(\alpha)\ne\emptyset$,
then it is smooth of dimension $\beta(n,d,k)$;
\item[(iii)]
if $U(\alpha)\ne\emptyset$, then it is irreducible; hence
$\overline{U(\alpha)}$ is an irreducible component of
$G(\alpha;n,d,k)$;
\item[(iv)] if $n\ge2$, then
$U(\alpha)\ne\emptyset$ if and only if the following conditions
hold: \be\label{eq:cond} 
(n-k)\alpha<d,\ k\le n+\frac1g(d-n)\ \mbox{and}\ (n,d,k)\ne(n,n,n);
\end{equation}
\item[(v)] if $n\ge2$, the set $U$ of Lemma
\ref{lemma:useful} is non-empty if and only if conditions
{\em(\ref{eq:cond2})} hold.
\end{itemize}
\end{theorem}
\pf
(i) is standard, while  (ii) follows at once from Lemma \ref{lemma:inj}.\\

For (iii), suppose first that $k<n$ and define
$$
V(\alpha)=\{(E,V)\in G(\alpha;n,d,k)\,|\,(E,V)\ \mbox{is
torsion-free}\}.
$$
Then $V(\alpha)$ is again open in $G(\alpha;n,d,k)$ and every
element of $V(\alpha)$ is given by an extension (\ref{eq:tf}) with
$h^0(F^*)=0$ by Lemma \ref{lemma:extclass}. It now follows,
exactly as in the proof of \cite[Theorem 4.3]{BGN} that
$V(\alpha)$ is irreducible if it is non-empty.

Recall now that every irreducible component of $G(\alpha;n,d,k)$
has dimension greater than or equal to
 $\beta(n,d,k)$ (see \cite[Corollary 3.6]{BGMN}).
To complete the proof of (iii), it is therefore sufficient to show
that \be\label{eq:tordim} \dim(U(\alpha)\setminus
V(\alpha))<\beta(n,d,k).
\end{equation}
The argument for this is the same as that of the corresponding
result for Brill-Noether loci \cite[Chapitre 3, Th\'eor\`eme
A.1]{M1}. For the convenience of the reader, we give an outline
here.

The points of $U(\alpha)\setminus V(\alpha)$ are represented by
extensions \be\label{eq:tor} 0\lra V\otimes\cO\lra E\lra F\oplus
T\lra 0,
\end{equation}
where $F$ is a vector bundle and $T$ is a torsion sheaf of length
$t\ge1$. The extension classes are defined by $k$-tuples
$$
e_1,\ldots,e_k\in\ \Ext^1(F\oplus T,\cO);
$$
by Lemmas \ref{lemma:dualcs} and \ref{lemma:extclass}, the
$\alpha$-stability of $(E,V)$ implies that $e_1,\ldots,e_k$ are
linearly independent and that $h^0(F^*)=0$. We now simply have to
estimate dimensions; for convenience, we will assume that the
support of $T$ consists of $t$ distinct points (the other cases
are handled similarly).

Since $(E,V)$ is $\alpha$-stable, the only automorphisms of
$(E,V)$ are scalar multiples of the identity \cite[Proposition
2.2]{BGMN}; it follows that the dimension of the component of
$U(\alpha)\setminus V(\alpha)$ consisting of the coherent systems
of the form (\ref{eq:tor}) for a fixed value of $t$ is
\begin{align*}
(n-k)^2(g-1)+1+t&+\dim \Gr(k,\Ext^1(F\oplus T,\cO))-\min\dim \Aut(F\oplus T)+1\\
&=\beta(n,d,k)+t-\min\dim \Aut(F\oplus T)+1.
\end{align*}
To establish (\ref{eq:tordim}), we therefore need to show that
$$
\dim \Aut(F\oplus T)\ge t+2.
$$
This is clear since $\dim \Aut T=t$ and $\dim \Hom(F,T)=(n-k)t\ge
t\ge1$.

This completes the proof of (iii) when $k<n$.
For the case $k=n$, see \cite[Theorem 5.6 and its proof]{BGMN}.\\

It remains to prove (iv) and (v). When $k<n$, the necessity of the
conditions (\ref{eq:cond}) and (\ref{eq:cond2}) has already been
proved in Lemma \ref{lemma:dualcs}, \cite[Lemma 4.1]{BGMN} and
Lemma \ref{lemma:extclass}. For $k=n$, (\ref{eq:cond}) and
(\ref{eq:cond2}) both reduce to $d>n$, which is a necessary
condition for non-emptiness by \cite[Remark 5.7]{BGMN}.

Finally, we shall prove that, if (\ref{eq:cond2}) holds, then
$U\ne\emptyset$. Since $U\subset U(\alpha)$ for all allowable
$\alpha$ by Proposition \ref{prop:restate} and Lemma
\ref{lemma:useful}, this will show also that
$U(\alpha)\ne\emptyset$ whenever (\ref{eq:cond}) holds.

For $k<n$, we note that, by Proposition \ref{prop:restate},
$G_L(n,d,k)\ne\emptyset$ and every element of it is torsion-free.
We claim that there exists a torsion-free coherent system $(E,V)$
for which $E$ is stable. Once this is proved, it follows that
$(E,V)$ arises from an extension (\ref{eq:tf}) with $h^0(F^*)=0$.
As already noted earlier in the proof, the set of all such
extensions is parametrised by an irreducible variety. It follows
from the openness of the stability condition that the general
extension (\ref{eq:tf}) defines a coherent system $(E,V)\in
G_L(n,d,k)$ with $E$ stable. Hence $U\ne\emptyset$.

For $d\le n$, the claim is proved in \cite{BGN}. For $d>n$, the
result can be deduced from a combination of \cite{BGN}, \cite{M1}
and \cite{BMNO}, and probably also from \cite{T}, but is perhaps
most easily obtained by using \cite{M3}; in fact it follows
directly from \cite[Th\'eor\`eme A.5]{M3}.

When $k=n$, the proof is on the same lines but using extensions of
the form (\ref{eq:tor1}) in place of those of the form
(\ref{eq:tf}). Again it is clear that the extensions are
parametrised by an irreducible variety. The fact that there exists
an extension (\ref{eq:tor1}) with $E$ stable is again a
consequence of \cite[Th\'eor\`eme A.5]{M3}. \qed

As a consequence of Theorem \ref{th:inj}, we obtain our first
important result on the geometry of $G(\alpha;n,d,k)$.

\begin{theorem}\label{th:alphai} Suppose $n\ge2$, $0<k\le n$ and
$\alpha>\alpha_I$. If the moduli space $G(\alpha;n,d,k)$ is
non-empty, then it is smooth and irreducible of dimension
$\beta(n,d,k)$ and is birationally equivalent to $G_L(n,d,k)$.
Moreover $G(\alpha;n,d,k)$ is non-empty if and only if the
conditions {\em(\ref{eq:cond})} hold.
\end{theorem}

\pf It follows from Definition \ref{def:inj} that, if
$\alpha>\alpha_I$, then $G(\alpha;n,d,k)=U(\alpha)$. Hence, if
$G(\alpha;n,d,k)\ne\emptyset$, Theorem \ref{th:inj}(v) implies
that $U\ne\emptyset$. The rest of the theorem now follows easily
from Theorem \ref{th:inj}. \qed

\begin{remark}\label{rmk3}
{\em If we denote by $G_I(n,d,k)$ the moduli space of coherent
systems of type $(n,d,k)$ which are $\alpha$-stable for $\alpha$
slightly greater than $\alpha_I$, the theorem can be restated to
say that $G_I(n,d,k)$ is birationally equivalent to $G_L(n,d,k)$.}
\end{remark}

\begin{corollary}\label{cor:birat}
Suppose $0<k< n$. If further $\alpha_I<\alpha<\frac{d}{n-k}$ and $k\le
n+\frac1g(d-n)$, then $G(\alpha;n,d,k)$ is birationally equivalent
to a fibration over the moduli space $M(n-k,d)$ with fibre the
Grassmannian $\Gr(k,d+(n-k)(g-1))$.
\end{corollary}

\pf This follows at once from the theorem and Proposition
\ref{prop:restate}. \qed

\begin{corollary}\label{cor:birat2}
Suppose $2\le k\le n$ and $d\le\min\{2n,n+\frac{ng}{k-1}\}$. Then
$G_0(n,d,k)$ is birationally equivalent to $G_L(n,d,k)$.
\end{corollary}

\pf This follows from Remark \ref{rmk3} and Corollary
\ref{cor:inj}.\qed

\section{The Diagram Lemma}\label{section:diag}

It follows from Theorem \ref{th:inj} that the ``flips'' at
critical points $\alpha>\alpha_I$ all have positive codimension.
The purpose of the next few sections is to obtain more information
about the flips when $k<n$ and $\alpha$ is large.

In this section we give a structural result that applies in
particular to all coherent systems which are $\alpha$-stable for
some $\alpha$ in the range
\begin{equation}\label{eq:range}
\max\left\{\frac{d-n}{n-k},0\right\}=\alpha_T < \alpha < \frac{d}{n-k}.
\end{equation}
By Lemma \ref{lemma:inj2}, all such coherent systems are
torsion-free. The result of this section may thus be viewed as an
extension of \cite[p.660, diagram (5)]{BGN} to $\alpha$-stable
coherent systems in the range (\ref{eq:range}) but without any
restriction on $d$.

\begin{lemma}\label{lemma:Main} {\bf (Diagram Lemma)}
Suppose that $k<n$ and that {\em(\ref{eq:range})} holds. Let $(E,V)$ be a
torsion-free
coherent system with $h^0(E^*)=0$. Suppose further that there
exists an exact sequence of coherent systems
 \begin{equation}\label{destab}
 0\to (E_1,V_1)\to (E,V)\to (E_2,V_2) \to 0
 \end{equation}
with $E_1$, $E_2$ both of positive rank, $h^0(E_1^*)=0$,
$(E_2,V_2)$ $\alpha$-semistable and $\mu_\alpha(E_2,V_2)\leq
\mu_\alpha (E,V)$. Then there exists a diagram

\be\label{array}
\begin{array}{ccccccccc}
&& 0 & & 0 & & 0 && \\
 &&  \downarrow & & \downarrow & & \downarrow && \\
0&\to & V_1\otimes{\cO} & \to & E_1 &\to & F_1&\to&0 \\
  && \downarrow & & \downarrow & & \downarrow &&\\
0&\to & V\otimes{\cO} & \to & E &\to & F &\to&0\\
  && \downarrow & & \downarrow & & \downarrow &&\\
0&\to &  V_2\otimes{\cO} & \to & E_2 &\to & F_2&\to&0 \\
   && \downarrow & & \downarrow & & \downarrow &&\\
&& 0 &  & 0 & & 0 && \\
\end{array}
\end{equation}

\noi where

\begin{itemize}
\item[(a)] the quotients $F_1$, $F$, and $F_2$ are all
locally free with positive rank,
\item[(b)]
$h^0(F_1^*)=h^0(F^*)=h^0(F_2^*)=0$,
\item[(c)] the extension
classes of $E_1$, $E$, $E_2$ are given respectively by $k_1$, $k$,
$k_2$ linearly independent vectors in $H^1(F_1^*)$, $H^1(F^*)$,
$H^1(F_2^*)$.
\end{itemize}

\end{lemma}

\pf Note first that, given $k<n$,

\begin{equation}\label{eq:slope}
\alpha<\frac{d}{n-k}\Leftrightarrow \alpha < \mu_\alpha (E,V)
\end{equation}
and
\begin{equation}\label{eq:slope2}
\frac{d-n}{n-k}<\alpha\Leftrightarrow\mu_\alpha (E,V)-1 < \alpha.
\end{equation}

The assumption that $\mu_\alpha (E_2,V_2) \leq \mu_\alpha (E,V)$
implies, by (\ref{eq:range}) and (\ref{eq:slope2}), that
$$
\mu_\alpha (E_2,V_2)-1\le\mu_\alpha(E,V)-1 < \alpha.
$$
Since $(E_2,V_2)$ is $\alpha$-semistable, this implies that
$(E_2,V_2)$ cannot possess a subsystem $(L,V_L)$ with $\deg
L\ge1$, $\rk L=1$ and $\dim V_L\ge1$. Now the proof of Lemma
\ref{lemma:inj2}(ii) shows that $(E_2,V_2)$ is torsion-free. On
the other hand, $(E_1,V_1)$ is a subsystem of the torsion-free
coherent system $(E,V)$ and is therefore torsion-free. An easy
application of the snake lemma
now gives us (\ref{array}).\\

If $F_1=0$, then $E_1\simeq\cO^{\oplus k_1}$ and
$\mu_{\alpha}(E_1,V_1)=\alpha$. But by hypothesis and
(\ref{eq:slope}),
$$
\mu_{\alpha}(E_1,V_1)\geq \mu_\alpha (E,V)>\alpha,
$$
giving a contradiction. Hence $F_1$ has positive rank.

If $F_2=0$, then $E_2\simeq\cO^{\oplus k_2}$, which contradicts the
assumption that $h^0(E^*)=0$.  This completes the proof of (a).

Finally we have $h^0(E^*)=h^0(E_1^*)=0$ by hypothesis and
$h^0(E_2^*)=0$ by considering the middle column of (\ref{array}).
Condition (b) follows by considering the rows of (\ref{array}),
and (c) by Lemma \ref{lemma:extclass}. \qed

\begin{remark}\label{rmk10}
{\em We have stated Lemma \ref{lemma:Main} under rather general
hypotheses. The situation in which we shall be applying it is that
of \cite[Lemma 6.5]{BGMN}. Slightly modifying the notation of
\cite{BGMN}, we suppose
\begin{itemize}
\item  $\alpha>\alpha_T$ is a critical value, \item $\alpha^+$
denotes a value of $\alpha$ slightly greater than $\alpha$, while
$\alpha^-$ denotes a value slightly less than $\alpha$, \item
$G^{\pm}(\alpha)$ is the set of points in $G(\alpha^{\pm};n,d,k)$
represented by coherent systems which are $\alpha^{\pm}$-stable
but not $\alpha^{\mp}$-stable.
\end{itemize}
If now $(E,V)$ is strictly $\alpha$-semistable and
$\alpha^+$-stable, we have (by \cite[Lemma 6.5]{BGMN}) an
extension (\ref{destab}) in which
\begin{itemize}
\item[(d)] $(E_1,V_1)$ and $(E_2,V_2)$ are $\alpha^+$-stable, with
$\mu_{\alpha^+}(E_1,V_1)<\mu_{\alpha^+}(E_2,V_2)$, \item[(e)]
$(E_1,V_1)$ and $(E_2,V_2)$ are $\alpha$-semistable, with
$\mu_{\alpha}(E_1,V_1)=\mu_{\alpha}(E_2,V_2)$.
\end{itemize}
If $(E,V)$ is strictly $\alpha$-semistable and $\alpha^-$-stable
then we have an extension of the same form, but with $\alpha^+$
replaced by $\alpha^-$ in the condition (d).

It follows from Lemma \ref{lemma:dualcs} that, in either case, all
the hypotheses of the Diagram Lemma are satisfied, so that we have
a diagram (\ref{array}) and conditions (a), (b), (c) hold as well
as (d) and (e).}
\end{remark}

\section{The case $k=n-1$}\label{section:n-1}

As a first application of the Diagram Lemma, we consider the case
$k=n-1$, where $n\ge2$. In this case $\alpha_T=\max\{d-n,0\}$ and
(\ref{eq:range}) becomes \be\label{eq:rangen-1}
\max\{d-n,0\}<\alpha<d.
\end{equation}

\begin{proposition}\label{prop:k=n-1}
Let $(E,V)$ be a coherent system of type $(n,d,n-1)$.  If $(E,V)$
is $\alpha$-stable for some $\alpha$ in the range
{\em(\ref{eq:rangen-1})}, then $(E,V)$ is $\alpha$-stable for all
$\alpha$ in the range {\em(\ref{eq:rangen-1})}.
\end{proposition}

\pf If the result is false, there is a critical value of $\alpha$
in the range (\ref{eq:rangen-1}) and a coherent system $(E,V)$
which is strictly $\alpha$-semistable. By \cite[Lemma 6.5]{BGMN},
there exists an extension (\ref{destab}) with $(E_1,V_1)$,
$(E_2,V_2)$ either both $\alpha^+$-stable or both
$\alpha^-$-stable and $\mu_\alpha(E_1,V_1)=\mu_\alpha(E_2,V_2)$.
It follows now from Lemma \ref{lemma:dualcs} that all the
hypotheses of the Diagram Lemma are satisfied. Applying this
lemma, we obtain $k_1\le n_1-1$ and $k_2\le n_2-1$. Thus
$$
k=k_1+k_2\le n_1+n_2-2=n-2\ ,
$$
which contradicts our assumption that $k=n-1$. \qed

\begin{corollary}\label{cor:k=n-1}
Suppose that {\em(\ref{eq:rangen-1})} holds. Then
$$
G(\alpha;n,d,n-1)=G_L(n,d,n-1).
$$
\end{corollary}

\pf This follows immediately from the proposition. \qed

\begin{theorem}\label{th:gt}
Suppose that {\em(\ref{eq:rangen-1})} holds. Then $G(\alpha;n,d,n-1)$
is non-empty if and only if $d\ge \max\{1,n-g\}$. When this
condition holds, $G(\alpha;n,d,n-1)$ is a fibration over the
Jacobian $J^d$, with fibre the Grassmannian $\Gr(n-1,d+g-1)$. In
particular
\begin{itemize}
\item[(i)] if $d>n-g$, then
$
\Pic(G(\alpha;n,d,n-1))\cong\Pic J^d\times \ZZ
$;
\item[(ii)] if $n>g$, then $\Pic(G(\alpha;n,n-g,n-1))\cong\Pic J^{n-g}$;
\item[(iii)] in all cases, $\Pic^0(G(\alpha;n,d,n-1))$ is isomorphic
to the Jacobian $J(C)$.
\end{itemize}
\end{theorem}

\pf (i) follows from Proposition \ref{prop:restate} and Corollary
\ref{cor:k=n-1}. For (ii), note that $G_L(n,n-g,n-1)\cong J^{n-g}$.
Finally, for (iii), recall that $J^d\cong J(C)$ and $\Pic^0(J(C))\cong J(C)$
because $J(C)$ is a principally polarized abelian variety. \qed

\begin{remark}\label{rmk21} {\em We have presented these results both in terms
of the {\em Picard group} $\Pic$ (the group of all algebraic line bundles)
and the
{\em Picard variety} $\Pic^0$ (the group of all
topologically trivial line bundles). Our reasons for doing this are firstly
that the results for Picard varieties are particularly simple
(see also Theorem \ref{picvar}) and
secondly that, in the coprime case, the Picard variety
has the structure of a polarized abelian variety and carries
a lot of geometrical information; the isomorphism of Theorem \ref{th:gt} is
certainly an isomorphism of abelian varieties. If one could prove that
it is an isomorphism of {\em polarized} abelian varieties, it would follow
that, under the conditions of the theorem, the variety $G(\alpha;n,d,n-1)$
satisfies a {\em global Torelli theorem}, that is, it determines $C$ as an
algebraic curve.}
\end{remark}

If we now denote by $G_T(n,d,k)$ the moduli space of
$\alpha$-stable coherent systems for $\alpha$ slightly greater
than $\alpha_T$, we have the following immediate corollary of Theorem
\ref{th:gt}.

\begin{corollary}\label{cor:gt}
If $d\ge \max\{1,n-g\}$, then $G_T(n,d,n-1)$ is a fibration over
the Jacobian $J^d$, with fibre the Grassmannian $\Gr(n-1,d+g-1)$.
\end{corollary}

\begin{remark}\label{rmk5}
{\em In fact, we can identify the fibration of Theorem \ref{th:gt}
and Corollary \ref{cor:gt} precisely. Let $d>0$ and let $\cP$ be a
Poincar\'e bundle on $C\times J^d$, that is a line bundle whose
restriction to $C\times\{j\}$ is the line bundle $L_j$ on $C$
corresponding to the point $j\in J^d$. Let $p_C$, $p_J$ be the
projections of $C\times J^d$ on its factors. Then the direct image
$p_{J*}\cP^*$ is zero and hence $R^1p_{J*}\cP^*$ is locally free
of rank $d+g-1$. The corresponding vector bundle classifies the
extensions
$$
0\lra\cO\lra E\lra L_j\lra0,
$$
where $j$ is a (variable) point of $J^d$. It follows that
$G_L(n,d,n-1)$ can be identified with the Grassmann bundle
$\Gr(n-1,R^1p_{J*}\cP^*)$ whose fibre over $j$ is the Grassmannian
of subspaces of $H^1(L_j^*)$ of dimension $n-1$. By relative Serre
duality, we can identify this with $\Gr(d+g-n,p_{J*}(\cP\otimes
p_C^*K_C))$, where $K_C$ is the canonical bundle. This variety is
well-known; it is the variety of linear systems of degree $d+2g-2$
and dimension $d+g-n-1$, classically denoted by
$G_{d+2g-2}^{d+g-n-1}$. The vector bundle $p_{J*}(\cP\otimes
p_C^*K_C)$ is an example of a {\em Picard bundle}; these are well
understood (for example, their Chern classes are known) and there
is a substantial literature on them (see, for instance \cite{K,
Ma1, Ma2, Mac, Mu}).}
\end{remark}

We finish this section with an application to Brill-Noether loci.

\begin{theorem}\label{th:bnn-1}
Let $d>0$ and $n-g\le d<n$. Then $B(n,d,n-1)$ can be identified
with $G_{d+2g-2}^{d+g-n-1}$.
\end{theorem}

\pf Since $d\le n$, $G_T(n,d,n-1)=G_0(n,d,n-1)$. By \cite[section 2.3]{BGMN}, we have a morphism
$$
\psi:G_0(n,d,n-1)\lra \widetilde{B}(n,d,n-1).
$$
By Remark \ref{rmk5}, the theorem will follow if we can prove that
$\psi$ is an isomorphism onto $B(n,d,n-1)$.

Note first that, by \cite[Theorems B and \~B]{BGN}, the stated conditions
imply that
$$
B(n,d,n-1)\ne\emptyset\ \mbox{and}\ \widetilde{B}(n,d,n)=\emptyset
$$
(for the second condition here, we need $d<n$). Now suppose that
$E$ is a strictly semistable bundle of rank $n$ and degree $d$
with $h^0(E)=n-1$. Then there exists a proper subbundle $E_1$ of
$E$ such that both $E_1$ and $E_2=E/E_1$ are semistable and have
the same slope. If we let $n_1$, $n_2$ denote the ranks of $E_1$,
$E_2$, then one of the $E_i$ must have $h^0(E_i)\ge n_i$. This
contradicts \cite[Theorem B]{BGN} since $\mu(E_i)=\mu(E)$ and
hence $0<\deg E_i<n_i$. So $B(n,d,n-1)=\widetilde{B}(n,d,n-1)$.
The result now follows from \cite[Corollary 11.5]{BGMN}. \qed

\section{Dimension counts and flips}\label{section:flip}

In this section we will again suppose that $k<n$ and obtain lower
bounds on the codimensions of the flips at all critical values
$\alpha>\alpha_T$. We use the notation of Remark \ref{rmk10}. We
are interested in the extensions (\ref{destab}) and note that, as
in \cite[equation (8)]{BGMN},
$$
\dim \Ext^1((E_2,V_2),(E_1,V_1))= C_{21}+\dim {\HH}^0_{21}+ \dim
{\HH}^2_{21},
$$
where
\begin{eqnarray}
{\HH}_{21}^0 &:= & \Hom((E_2,V_2),(E_1,V_1))\ ,\nonumber\\
\HH^2_{21}&:= & \Ext^2((E_2,V_2),(E_1,V_1))\ ,\nonumber\\
C_{21}&:= & k_2\chi (E_1)-\chi(E_2^{*}\otimes E_1)- k_1k_2
\nonumber\\
& = &  n_1n_2(g-1)-d_1n_2+d_2n_1+k_2d_1-k_2n_1(g-1)-k_1k_2\ .
\label{C21}
\end{eqnarray}

\begin{remark}\label{rmk2}
 {\em Notice that, if
$\mu_{\alpha^+}(E_1,V_1)<\mu_{\alpha^+}(E_2,V_2)$
 and $\mu_{\alpha}(E_1,V_1)=\mu_{\alpha}(E_2,V_2)$, then it follows that
$\mu_{\alpha^-}(E_1,V_1)>\mu_{\alpha^-}(E_2,V_2)$. In this case
$(E_1,V_1)$ is an $\alpha^-$-destabilizing subsystem, or,
equivalently,  $(E_2,V_2)$ is an $\alpha^-$-destabilizing quotient
system. Similarly, $(E_1,V_1)$ and $(E_2,V_2)$ are
$\alpha^+$-destabilizing for $(E,V)$ if
$\mu_{\alpha^-}(E_1,V_1)<\mu_{\alpha^-}(E_2,V_2)$.}
\end{remark}

\begin{lemma}\label{lemma:h0=h2=0}
Suppose that we have a diagram {\em(\ref{array})} and that condition {\em(d)}
in Remark \ref{rmk10} holds.
Then ${\HH}_{21}^0 = 0 = \HH^2_{21}$ and thus
\begin{equation}\label{dim-ext}
\dim \Ext^1((E_2,V_2),(E_1,V_1)) = C_{21}.
\end{equation}
\end{lemma}

\pf ${\HH}_{21}^0 = 0$ follows at once from condition (d).
Moreover ${\HH}_{21}^2 = 0$ because (see
 \cite[Proposition 3.2]{BGMN})
 $$
 {\HH}_{21}^2 = H^0(E_1^*\otimes N_2\otimes K)^*,
 $$
\noi where $N_2$ is the kernel of the map
$V_2\otimes\cO\longrightarrow
 E_2$. But $N_2=0$ since $(E_2,V_2)$ is injective.

\qed

By Lemma \ref{lemma:inj}, $G(\alpha^{\pm};n_i,d_i,k_i)$ have their
expected dimensions. Together with Lemma \ref{lemma:h0=h2=0} this
yields lower bounds (see \cite[section 6.3]{BGMN})

\be\label{eqtn:codim} \codim
G^{\pm}(\alpha)\geq\min\left\{C_{12}\right\}\ ,
\end{equation}

\noi where

\begin{equation}\label{C12}
C_{12}:=  n_2n_1(g-1)-d_2n_1+d_1n_2+k_1d_2-k_1n_2(g-1)-k_2k_1
\end{equation}

\noi and the minimum is taken over all possible extension types
which can occur for coherent systems in $G^{+}(\alpha)$
(respectively $G^{-}(\alpha)$).

\begin{lemma}\label{lemma:constraints} For $G^{+}(\alpha)$,
the minimum in {\em(\ref{eqtn:codim})} is taken over all
$(n_1,d_1,k_1)$ and $(n_2,d_2,k_2)$ satisfying
\begin{enumerate}
\item[{\rm(i)}] $n_1+n_2=n,\ d_1+d_2=d,\ k_1+k_2=k$,
\item[{\rm(ii)}] $0<n_1<n$, \item[{\rm(iii)}] $0\le
\frac{k_1}{n_1}<\frac{k}{n}<\frac{k_2}{n_2}<1$, \item[{\rm(iv)}]
${\displaystyle\alpha=\frac{nd_1-n_1d}{n_1k-nk_1}}$.
\end{enumerate}

\noi Similarly, for $G^{-}(\alpha)$, the minimum is taken over all
$(n_1,d_1,k_1)$ and $(n_2,d_2,k_2)$ satisfying
\begin{enumerate}
\item[{\rm(i)}$'$] $n_1+n_2=n,\ d_1+d_2=d,\ k_1+k_2=k$,
\item[{\rm(ii)}$'$] $0<n_2<n$, \item[{\rm(iii)}$'$] $0\le
\frac{k_2}{n_2}<\frac{k}{n}<\frac{k_1}{n_1}<1$,
\item[{\rm(iv)}$'$] $\displaystyle{\alpha=\frac{nd_2-n_2d}{n_2k-nk_2}}$.
\end{enumerate}
\end{lemma}

\pf We prove only the first set of constraints, i.~e. (i)--(iv).
The proof for the other set is similar.

Constraints (i) and (ii) are obvious and (iv) is a restatement of
the condition $\mu_{\alpha}(E_1,V_1)=\mu_\alpha(E,V)$. Since
coherent systems in $G^{+}(\alpha)$ are $\alpha$-semistable but
$\alpha^{+}$-stable, constraint (iii)  follows from the Diagram
Lemma. \qed

\begin{proposition}\label{prop:goodflip} Let $\alpha$ be a critical
value in the range
$$\alpha_T<\alpha<\frac{d}{n-k}\ .$$
\noi Then

\be\label{eq:C12-} C_{12}\ge(n_1-k_1)(n_2-k_2)(g-1)+1\ge g
\end{equation}

\noi for all extension types (as in Lemma \ref{lemma:Main}) which
can occur for coherent systems in $G^{-}(\alpha)$, and

\be\label{eq:C12+}
C_{12}\ge(g-1)(n_1-k_1)(n_2-k_2)+d_1n_2-d_2n_1+1\ge g+1
\end{equation}

\noi for all extension types (as in Lemma \ref{lemma:Main}) which
can occur for coherent systems in $G^{+}(\alpha)$.
\end{proposition}

\pf Note first that, since the coherent systems $(E_i,V_i)$
satisfy the hypotheses of Lemma \ref{lemma:extclass} with
appropriate choice of either $\alpha^+$ or $\alpha^-$, we have
(restating (\ref{eq:kbound})) \be\label{fact2}
d_i-k_i\ge(k_i-n_i)(g-1).
\end{equation}

\subsubsection*{Case 1: $G^-(\alpha)$}

The equation $\mu_\alpha(E_1,V_1)=\mu_\alpha(E_2,V_2)$ can be
rewritten as \be\label{eq:12}
\left(\frac{d_1}{n_1-k_1}-\alpha\right)=
\frac{n_1(n_2-k_2)}{n_2(n_1-k_1)}\left(\frac{d_2}{n_2-k_2}-\alpha\right).
\end{equation}
Note also that
$$
\alpha<\frac{d}{n-k}=\frac{d_1+d_2}{(n_1-k_1)+(n_2-k_2)}.
$$
It follows that at least one side of (\ref{eq:12}) is positive.
Hence both sides are positive. Since also
$n_1(n_2-k_2)>n_2(n_1-k_1)>0$ by (iii)$'$, it follows from
(\ref{eq:12}) that
$$
 \frac{d_1}{n_1-k_1}>\frac{d_2}{n_2-k_2}\ ,
$$
i.~e.
$$
d_1n_2-d_2n_1+k_1d_2>k_2d_1\ .
$$
(\ref{C12}) and (\ref{fact2}) for $i=1$ now give
$$
C_{12}> n_2(n_1-k_1)(g-1)+k_2(d_1-k_1)\ge(n_1-k_1)(n_2-k_2)(g-1) .
$$
This gives (\ref{eq:C12-}).

\subsubsection*{Case 2: $G^+(\alpha)$}

Using (\ref{C12}) and (\ref{fact2}) for $i=2$, we get

\begin{eqnarray*}
C_{12}& \ge& (g-1)(n_2 (n_1-k_1) +k_1(k_2-n_2))
+d_1n_2-d_2n_1\\
&=& (g-1)(n_1n_2-2n_2k_1+k_1k_2)+d_1n_2-d_2n_1\ .
\end{eqnarray*}

\noi But, since $n_1k_2>n_2k_1$ by (iii), we get
$(n_1n_2-2n_2k_1+k_1k_2)>(n_1-k_1)(n_2-k_2)$, and hence

$$
C_{12} > (g-1)(n_1-k_1)(n_2-k_2)+d_1n_2-d_2n_1.
$$
This gives the first inequality in (\ref{eq:C12+}). For the
second, note that $\frac{k_2}{n_2}>\frac{k_1}{n_1}$ and
$\mu_\alpha(E_1,V_1)=\mu_\alpha(E_2,V_2)$ give
$\frac{d_1}{n_1}>\frac{d_2}{n_2}$. \qed

We know already from Theorem \ref{th:alphai} that all flips in the
range $\alpha>\alpha_I$ are good. Proposition \ref{prop:goodflip}
applies to a more restricted range for $\alpha$ but gives a much
stronger conclusion.
\begin{theorem}\label{codim-g}
Suppose that $0<k<n$ and that {\em(\ref{eq:range})} holds. Then,
$G(\alpha;n,d,k)$ and $G_L(n,d,k)$ are smooth and irreducible. Moreover
these varieties are isomorphic outside subvarieties  of codimension
at least $g$.
\end{theorem}
\pf Smoothness and irreducibility have already been proved
(Theorem \ref{th:alphai}).
The rest follows from Proposition \ref{prop:goodflip}. \qed

\section{Picard Groups and  homotopy groups for $k<n$}\label{section:picard}

We begin with the following key proposition.

\begin{proposition}\label{picard-homotopy}
Let $0<k<n$ and $d>0$. Suppose further that $k\leq n+\frac1g(d-n)$.
Then, for $\alpha_T<\alpha<\frac{d}{n-k}$,

\begin{itemize}
\item[(i)] $ \Pic(G(\alpha;n,d,k))\cong \Pic(G_L(n,d,k))$;
\item[(ii)] $\pi_i(G(\alpha;n,d,k))\cong \pi_i(G_L(n,d,k)) \;\;\;
\mbox{for}\;\;\; i\leq 2g-2$.
\end{itemize}
\end{proposition}
\pf This follows from Theorem \ref{codim-g}. \qed

In order to apply this, we need to calculate the Picard groups and
homotopy groups of $G_L$. We suppose first that $n-k$ and $d$ are coprime.

\begin{theorem}\label{picard}
Let $0<k<n$ and $d>0$. Suppose that $n-k$ and $d$ are coprime and
that $k\le n+\frac1g(d-n)$. Then
\begin{itemize}
\item[(i)] if $k< n+\frac1g(d-n)$, then
$ \Pic(G_L(n,d,k)) \cong  \Pic(M(n-k,d))\times \ZZ$;
\item[(ii)] if $k= n+\frac1g(d-n)$, then $\Pic(G_L(n,d,k)) \cong  \Pic(M(n-
k,d))$.
\end{itemize}
\end{theorem}
\pf (i) follows from the fact that, by Proposition
\ref{prop:restate}(iv),
 $G_L$ is a fibration over $M(n-k,d)$ with fibre $\Gr(k,d+(n-k)(g-1))$.
(Note that the assumption on $k$ implies that this Grassmannian has
positive dimension.)

(ii) is clear since, in this case $G_L(n,d,k) \cong  M(n-k,d)$.
\qed

Combining  Theorem \ref{picard} with Proposition
\ref{picard-homotopy}(i) we have the following.

\begin{corollary}\label{cor:picard}
Let $0<k<n$ and $d>0$. Suppose that $n-k$ and $d$ are coprime and
that $k\le n+\frac1g(d-n)$ and $\alpha_T<\alpha<\frac{d}{n-k}$. Then
\begin{itemize}
\item[(i)] if $k< n+\frac1g(d-n)$, then
$ \Pic(G(\alpha;n,d,k)) \cong  \Pic(M(n-k,d))\times \ZZ$;
\item[(ii)] if $k= n+\frac1g(d-n)$, then $\Pic(G(\alpha;n,d,k)) \cong
\Pic(M(n-k,d))$.
\end{itemize}
\end{corollary}

To compute the homotopy groups of $G_L=G_L(n,d,k)$,
one can  use the
homotopy sequence for $G_L$ as a ${\Gr}(k,N)$ fibration over
$M(n-k,d)$  with  $N=d+(n-k)(g-1)$, given by

\begin{equation}\label{homotopy-sequence}
 ...  \lra  \pi_i(\Gr(k,N)) \lra  \pi_i(G_L)  \lra  \pi_i(M(n-k,d)) \lra
\pi_{i-1}(\Gr(k,N)) \lra ...
\end{equation}

\begin{theorem}\label{homotopy:k=n-1}
Let   $d\ge\max\{1,n-g\}$ and suppose that $\max\{d-n,0\}<\alpha <d$. Then

\begin{itemize}
\item[(i)] $\pi_i (G(\alpha;n,d,n-1)) \cong
\pi_i(\Gr(n-1,d+g-1))$ for $i\ge0$, $i\neq 1$;

\item[(ii)] $\pi_1 (G(\alpha;n,d,n-1)) \cong  H_1(C, \ZZ)$.
\end{itemize}
\end{theorem}

\pf By  Theorem \ref{th:gt}, $G_L$ is a $\Gr(n-1,d+g-1)$ fibration
over $J^d$. The result follows now from (\ref{homotopy-sequence})
(with $k=n-1$)
and  the fact that $\pi_i(J^d)= 0$ for $i\neq 1$ and
$\pi_1(J_d)=H_1(C;\ZZ)$. \qed

\begin{corollary}
Let $d>0$ and $n-g\leq d \leq n$. Then, for $0<\alpha <d$,

\begin{itemize}

\item[(i)] $\pi_i (G(\alpha;n,d,n-1)) \cong
\pi_i(\Gr(n-1,d+g-1))$ for $i\ge0$, $i\neq 1$;

\item[(ii)] $\pi_1 (G(\alpha;n,d,n-1))\cong H_1(C, \ZZ)$.
\end{itemize}
\end{corollary}

From the above results and \cite{BGN} one can derive for $n-g \leq
d<n$ the Picard group and homotopy groups  of the Brill-Noether locus
$B(n,d,n-1)$.
However, this also follows from the explicit description of
$B(n,d,n-1)$ given by Theorem \ref{th:bnn-1}. Note also that, if $d<n-g$,
then $G(\alpha;n,d,n-1)=\emptyset$ for all $\alpha$ by Lemma \ref{lemma:inj2}
and Theorem \ref{th:alphai}.

Building upon the gauge-theoretic  approach of Atiyah and Bott \cite{AB} to
the moduli space of stable bundles,  Daskalopoulos and Uhlenbeck
\cite{D,DU} have computed some of the homotopy groups of $M(r,d)$.
This information can be used to compute $\pi_1(G_L)$ and
$\pi_2(G_L)$.

\begin{theorem}\label{homotopy}
Let $0<k\le n-2$ and $d>0$. Suppose further that $n-k$ and $d$ are coprime
and that $\alpha_T<\alpha<\frac{d}{n-k}$. Then
\begin{itemize}

\item[(i)] $\pi_1(G(\alpha;n,d,k))\cong
\pi_1(M(n-k,d))\cong H_1(C,\ZZ)$;

\item[(ii)] if $k< n+\frac1g(d-n)$,
$\pi_2(G(\alpha;n,d,k))\cong  \ZZ \times \ZZ$;
\item[(iii)] if $k= n+\frac1g(d-n)$, $\pi_2(G(\alpha;n,d,k))\cong  \ZZ$.
\end{itemize}
\end{theorem}
\pf From  Proposition \ref{prop:restate}(iv), $G_L$ is a Grassmann
fibration over $M(n-k,d)$ and we have the homotopy sequence
(\ref{homotopy-sequence}). From this, since $\pi_0(\Gr(k,N))
=\pi_1(\Gr(k,N))=0$,
 we deduce that $\pi_1(G_L)\cong \pi_1(M(n-k,d))$.
It follows from \cite[Theorem 9.12]{AB} that $\pi_1(M(n-k,d))$
is isomorphic to $\pi_1(J^d)$, which
in turn is isomorphic to
$H_1(C,\ZZ)$. Statement (i) now follows from
Proposition \ref{picard-homotopy}(ii).

To compute the second homotopy group, we first note that
the condition $k< n+\frac1g(d-n)$ (resp.~$k= n+\frac1g(d-n)$) is
equivalent to $k<N$ (resp.~$k=N$). Moreover,
since
$\pi_1(\Gr(k,N))=0$, (\ref{homotopy-sequence}) gives
\begin{equation}\label{2-homotopy}
 ...  \lra  \pi_2(\Gr(k,N)) \lra  \pi_2(G_L)  \lra  \pi_2(M(n-k,d)) \lra 0.
\end{equation}

For $k<N$, we use again the fact that $\pi_1(\Gr(k,N))=0$ to deduce that
$$
\pi_2(\Gr(k,N))\cong H_2(\Gr(k,N),\ZZ)\cong \ZZ.
$$
We claim that the map
$f: \ZZ  \lra  \pi_2(G_L)$,
induced by these isomorphisms and (\ref{2-homotopy}), is injective.
This is true because $\Gr(k,N)$ is a
subvariety of $G_L$ and the map $H_2(\Gr(k,N),\ZZ)\lra
H_2(G_L,\ZZ)$ must be  injective (since the restriction of an ample
line bundle over $G_L$  to $\Gr(k,N)$  must give an ample line
bundle) and factors  through $f$. Now, from \cite{DU}, we have that
$\pi_2(M(n-k,d)) \cong\ZZ$, and from (\ref{2-homotopy}) we deduce that
$\pi_2(G_L)\cong\ZZ\times\ZZ$. Statement (ii) follows now from Proposition
\ref{picard-homotopy}(ii).

Finally, if $k=N$, $\Gr(k,N)$ is a point; (iii) now follows from
(\ref{2-homotopy}) and the fact that $\pi_2(M(n-k,d)) \cong\ZZ$.
\qed

\begin{remark}
{\em More information on higher homotopy groups could be obtained
from the knowledge of the higher homotopy groups of the moduli
space of stable bundles, which in turn are related, as shown in
\cite{DU}, to the homotopy groups of the unitary gauge group.}
\end{remark}

\begin{remark}{\em
A direct approach, similar to the one in \cite{D,DU}, to compute
the homotopy groups of the moduli space of coherent systems should
also be possible in general. In fact, this has been carried out
by Bradlow and Daskalopoulos \cite{BD}  for $k=1$, with some
additional restrictions on $d$ and $\alpha$, but no restrictions
on the coprimality of $n-1$ and $d$. They compute the first and
second homotopy groups, which coincide with the results given by
Theorem \ref{homotopy}.}
\end{remark}

We show  now how the results on the Picard group and
the first and second homotopy groups  are also valid, under certain
restrictions,
when $n-k$ and $d$ are not coprime. The principal assertion of
Proposition \ref{prop:restate}(iii) is that every coherent system
$(E,V)\in G_L$  defines an extension (\ref{eq:tf}) with $F$
semistable and $h^0(F^*)=0$. Hence, the key point is to find good estimates for
the
number of parameters counting strictly semistable bundles. The results of
\cite{DU} then allow us to carry out the necessary
computations.

To obtain our estimates we can use the Jordan-H\"older filtration
of a semistable bundle $F$, given by
\begin{equation}\label{j-h-filtration}
0=F_0\subset F_1\subset F_2 \subset ... \subset F_r=F,
\end{equation}
with $Q_i=F_i/F_{i-1}$ stable and
$\mu(Q_i)=\frac{d_i}{m_i}=\mu(F)$ for $1\leq i \leq r$. The bundle
$Q=\oplus_i Q_i$ is the {\em graduation} of $F$.  (Similar computations
to the ones given below can be found in \cite{AB,BGN,Ty1}.)

\begin{proposition} Let $\cS$ be a set of isomorphism classes of
semistable bundles of rank $m$ and degree $d$, whose Jordan-H\"older
filtration {\em(\ref{j-h-filtration})} has graduation $Q=\oplus_{i=1}^r Q_i$,
with $m_i=\rk Q_i$ and $d_i=\deg Q_i$. Then $\cS$ depends on at most
\begin{equation}\label{dim-family}
\left(\sum_i m_i^2+\sum_{i<j} m_im_j\right)(g-1) +1
\end{equation}
parameters.
\end{proposition}
\pf This is obtained by adding up the numbers
$$
\dim M(m_i,d_i)=m_i^2(g-1) +1
$$
for $1\leq i\leq r$, and the dimensions of the spaces of equivalence
classes of extensions
$$
0\lra F_{j-1}\lra F_j\lra Q_j \lra 0
$$
for $2\leq j\leq r$. By Riemann-Roch and the condition
$\mu(Q_i)=\mu(Q_j)$,  these dimensions are given by
$$
h^1(Q_j^*\otimes
F_{j-1})-1=m_j\left(\sum_{i<j}m_i\right)(g-1)-1.
$$
(We have  assumed that $Q_i$ and $Q_j$
are not isomorphic, since if they are isomorphic
the number of parameters on which $\cS$ depends is actually smaller.) \qed

\begin{corollary}\label{codim} Let $0<k<n$ and suppose that
$G_L(n,d,k)\ne\emptyset$. Let $W$ denote the closed subvariety of $G_L(n,d,k)$
consisting of coherent systems which arise from extensions {\em(\ref{eq:tf})}
in which $F$ is strictly semistable. Then
the codimension  of $W$ in $G_L(n,d,k)$ is at least
\begin{equation}\label{min}
 D:=\min \left\{\left(\sum_{i<j} m_im_j\right)(g-1) \right\},
\end{equation}
where the minimum is taken over all sequences of positive integers
$r,m_1,\ldots,m_r$ such that $r\ge2$ and $\sum m_i=n-k$.
\end{corollary}
\pf For fixed $F\in\cS$ with $h^0(F^*)=0$, the dimension of the variety
of coherent systems of the form (\ref{eq:tf}) is
$$
k\cdot h^1(F^*)-k^2=k(d+(n-k)(g-1)-k).
$$
Adding this to (\ref{dim-family}) and subtracting the total from
$\dim G_L(n,d,k)=\beta(n,d,k)$ gives the estimate (\ref{min}).
\qed
\begin{proposition} \label{prop-min}
The minimum in {\em(\ref{min})} is attained when $r=2$ and
$\{m_1,m_2\}=\{1,n-k-1\}$. Hence
$$
D=(n-k-1)(g-1).
$$
\end{proposition}
\pf The minimum of the expression $\sum_{i<j} m_im_j$ is achieved for
$r=2$. So
$$
D=\min \{ m_1(n-k-m_1)(g-1) \} = (n-k-1)(g-1).
$$ \qed

\begin{theorem}\label{picard-general}
Let $0<k<n$ and $d>0$. Suppose that $k<n+\frac1g(d-n)$ and
\begin{equation} \label{g-inequality}
(n-k-1)(g-1)\geq 2
\end{equation}
Then, for $\alpha_T<\alpha<\frac{d}{n-k}$,

\begin{itemize}

\item[(i)] $ \Pic(G(\alpha;n,d,k))\cong \Pic(M(n-k,d))\times
\ZZ$;

\item[(ii)] $\pi_1(G(\alpha;n,d,k))\cong
\pi_1(M(n-k,d))\cong H_1(C,\ZZ)$;

\item[(iii)] there is an exact sequence
$$0\lra\ZZ\lra\pi_2(G(\alpha;n,d,k))\lra \ZZ \times\ZZ_p\lra0,$$ 
where $p=\GCD(n-k,d)$.
\end{itemize}
\end{theorem}
\pf By Proposition \ref{prop:restate}(iii), $G_L\setminus W$ is isomorphic
to a Grassmann fibration over $M(n-k,d)$. It follows from
Proposition \ref{prop-min} and (\ref{g-inequality}) that
$$
\Pic(G_L)\cong\Pic(G_L\setminus W) \cong  \Pic(M(n-k,d))\times \ZZ.
$$
The statement (i) now follows from Proposition \ref{picard-homotopy}.

To prove
(ii) and (iii) we use the same arguments as in Theorem
\ref{homotopy} taking into account now that, as proved
in \cite[Theorem 3.1]{DU},
$$\pi_1(M(n-k,d))\cong H_1(C,\ZZ),\quad\pi_2(M(n-k,d))\cong  \ZZ \times \ZZ_p.
$$
This gives $\pi_1(G_L\setminus W)\cong\pi_1(M(n-k,d))$ and an exact sequence
$$0\lra\ZZ\lra\pi_2(G_L\setminus W)\lra \ZZ \times\ZZ_p\lra0$$
(compare (\ref{2-homotopy})). 
Now use Proposition \ref{prop-min} and (\ref{g-inequality}) again.
\qed

\begin{remark}\label{rmk22}
{\em When $k=n+\frac1g(d-n)$, $\Gr(k,N)$ is a point; in this case
(ii) remains true, but we must delete the factor $\ZZ$ in (i) and (iii) becomes
$\pi_2(G(\alpha;n,d,k))\cong \ZZ \times\ZZ_p$. It is a plausible conjecture, 
compatible with Theorem \ref{homotopy}, \cite[Theorem 1.12]{BD}, Theorem \ref{picard-general} and 
the first part of this remark, that, for $0<k\le n-2$, $\pi_2(G(\alpha;n,d,k))\cong\ZZ\times \ZZ \times\ZZ_q$, where
$q=\GCD(n,d,k)$, with one factor $\ZZ$ being dropped when $k=n+\frac1g(d-n)$.}
\end{remark}
\begin{remark}\label{rmk23}
{\em The inequality (\ref{g-inequality}) fails only for $k=n-1$ (any genus)
and for $k=n-2$ when $g=2$. The cases $k=n-1$ and $g=2$, $k=n-2$, $d$ odd
are covered by Theorems \ref{picard}, \ref{homotopy:k=n-1} and \ref{homotopy}; note that, if $k=n-1$, 
(i) and (ii) are true, but (iii) must be replaced by $\pi_2(G(\alpha;n,d,k))\cong \ZZ$. This leaves only the
case
$$
g=2,\ k=n-2,\ d\ \mbox{even},
$$ for which the theorem may fail. One may note that \cite{DU} specifically
excludes this case.}
\end{remark}

We finish with the result for the Picard variety, which takes a very simple
form.
\begin{theorem}\label{picvar} Let $0<k<n$ and $d>0$. Suppose that
$k\le n+\frac1g(d-n)$ and that $\alpha_T<\alpha<\frac{d}{n-k}$. Then,
except possibly in the case $g=2$, $k=n-2$, $d$ even,
$$
\Pic^0(G(\alpha;n,d,k))\cong J(C).
$$
\end{theorem}
\pf From \cite[Th\'eor\`emes A, C]{DN} and the fact that the
codimension of the strictly semistable locus in $M(n-k,d)$ is at least $2$,
it follows  that $\Pic^0(M(n-k,d))\cong J(C)$. For the rest,
see the proofs of Theorems \ref{picard} and \ref{picard-general}.
\qed

\section{The case $k=n-2$: Poincar\'e polynomials}\label{section:n-2}

So far, in applying the Diagram Lemma, we have used only the
numerical consequences of the construction. In this section we
make a first application of the geometry of the flips at critical
values in the range $\alpha>\alpha_T$. We shall be mainly concerned
with the case $k=n-2$; for values of $k<n-2$, some information can
be obtained from the Diagram Lemma, but additional techniques will
be necessary to obtain complete results.

We begin, however,  with a basic observation which applies to any  $k<n$
with $n-k$ and $d$ coprime. We write $P(X)$ for the Poincar\'e
polynomial of a space $X$ (with coefficients in any fixed field).

\begin{proposition}\label{poincare}
Let $0<k<n$ and $d>0$. Suppose that   $k\leq n+\frac1g(d-n)$ and that
$n-k$ and $d$ are coprime. Then
the Poincar\'e polynomial of $G_L$ is given by
$$
P(G_L(n,d,k))=P(M(n-k,d))\cdot P({\Gr}(k,N)),
$$
where $N=d+(n-k)(g-1)$.
\end{proposition}

\pf By Proposition \ref{prop:restate}(iv), $G_L(n,d,k)$ is the Grassmann
fibration over $M(n-k,d)$ associated with some vector bundle.
The result is now standard.
 \qed

\begin{remark} {\em (i) Note that, when $\GCD(n-k,d)=1$, $H^*(M(n-k,d);\ZZ)$
is torsion-free \cite[Theorem 9.9]{AB}. Hence $P(M(n-k,d))$ and
$P(G_L(n,d,k))$ are independent of the characteristic of the coefficient field.
A closed (though, in general, complicated) formula for $P(M(n-k,d))$ is
known \cite{Z}. Here we shall need only the special case
\begin{equation}\label{moduli}
P(M(2,d))(t)=\frac{(1+t)^{2g}((1+t^3)^{2g}-t^{2g}(1+t)^{2g})}{(1-t^2)(1-t^4)}
\ \ (d \mbox{ odd}),\end{equation}
which is essentially proved in \cite{N}.

(ii) For $0<k\le N$, the Poincar\'e polynomial of $\Gr(k,N)$ is
given by
\begin{equation}\label{grass}
P(\Gr(k,N))(t)=\frac{(1-t^{2(N-k+1)})(1-t^{2(N-k+2)})\cdots(1-t^{2N})}
{(1-t^2)(1-t^4)\cdots(1-t^{2k})}.
\end{equation}
Combining this with (i), one can obtain an explicit formula for
$P(G_L(n,d,k))$.}
\end{remark}

We suppose now that $k=n-2$. We are concerned with
critical values in the range (\ref{eq:range}), which in this case
is \be\label{eq:rangen-2} \max\{d-n,0\}<2\alpha<d.
\end{equation}
For convenience, we will write $G(\alpha)$, $G_L$ and $G_T$ for
$G(\alpha;n,d,n-2)$, $G_L(n,d,n-2)$ and $G_T(n,d,n-2)$.

At any such critical value $\alpha$, we know, either by Lemma
\ref{lemma:constraints} or
directly from the Diagram Lemma, that $k_1<n_1$ and $k_2<n_2$. Hence
\be\label{eq:kn-2} k_1=n_1-1,\ \ k_2=n_2-1.
\end{equation}
Inserting these into Lemma \ref{lemma:constraints} for
$G^+(\alpha)$, we get \be\label{eq:n1} 2n_1<n
\end{equation}
and \be\label{eq:alphan-2} \alpha=\frac{nd_1-n_1d}{n-2n_1}.
\end{equation}
Equation (\ref{eq:rangen-2}) now gives
$$
\max\{(d-n)(n-2n_1),0\}<2(nd_1-n_1d)<d(n-2n_1),
$$
i.e. \be\label{eq:d1}
\max\left\{d+2n_1-n,\frac{2n_1d}n\right\}<2d_1<d.
\end{equation}

\begin{lemma}\label{disjoint}
{\rm(i)} The flip locus $G^+(\alpha)$ is a disjoint union
$$
G^+(\alpha)=\bigsqcup G^+(n_1,d_1),
$$
where $(n_1,d_1)$ ranges over the set of possible solutions of
{\em(\ref{eq:n1})}, {\em(\ref{eq:alphan-2})}, {\em(\ref{eq:d1})}. Here
$G^+(n_1,d_1)$ is a smooth subvariety of $G(\alpha^+)$ and is
isomorphic to a projective bundle over $G_L(n_1,d_1,n_1-1)\times
G_L(n_2,d_2,n_2-1)$ with fibre dimension $C_{21}(n_1,d_1)-1$,
where \be\label{C211}
C_{21}(n_1,d_1)=n_1(g-1)+dn_1-d_1(n_1+1)-(n_1-1)(n-n_1-1).
\end{equation}
Moreover $G^+(n_1,d_1)$ has
codimension $C_{12}(n_1,d_1)$ in $G(\alpha^+)$ given  by
\be\label{C121}
C_{12}(n_1,d_1)=(n-n_1)(g-1)-d+d_1(n-n_1+1)-(n_1-1)(n-n_1-1).
\end{equation}

{\rm(ii)} Similarly
$G^-(\alpha)$ is a disjoint union
$$
G^-(\alpha)=\bigsqcup G^-(n_1,d_1)
$$
of smooth subvarieties of $G(\alpha^-)$, where $G^-(n_1,d_1)$ is
isomorphic to a projective bundle over $G_L(n_1,d_1,n_1-1)\times
G_L(n_2,d_2,n_2-1)$ with fibre dimension $C_{12}(n_1,d_1)-1$. Moreover
$G^-(n_1,d_1)$ has
codimension $C_{21}(n_1,d_1)$ in $G(\alpha^-)$.
\end{lemma}

\pf
For fixed $n$, $d$, the inequalities (\ref{eq:n1}) and
(\ref{eq:d1}) give rise to a finite number of choices for $n_1$,
$d_1$, hence also for $\alpha$. One can check that the values of
$\alpha$ lie in the torsion-free range for coherent systems of
types  $(n_1,d_1,n_1-1)$
and $(n_2,d_2,n_2-1)$. It follows from Proposition
\ref{prop:k=n-1} that $\alpha$ is not critical for either
of these types of coherent systems; hence the filtration
(\ref{destab}) of the Diagram Lemma is the unique Jordan-H\"older
filtration of $(E,V)$. It is possible that, for a given $\alpha$,
there may be more than one set of values for $n_1$, $d_1$
satisfying (\ref{eq:n1}), (\ref{eq:alphan-2}) and (\ref{eq:d1}).
However the uniqueness of the Jordan-H\"older filtration implies
that the flip locus $G^+(\alpha)$ is a disjoint union of subvarieties
as required. The smoothness of these subvarieties is proved as in
\cite[section 3]{Th} (see Appendix for details). 
The formulae (\ref{C211}) and (\ref{C121}) are obtained
from  (\ref{C21}) and (\ref{C12}) using (\ref{eq:kn-2})
and Lemma \ref{lemma:constraints}(i)).

(ii) is proved by interchanging the subscripts $12$ in the proof of (i)
\qed

This situation is analogous to that of Thaddeus \cite{Th}. Using
the same method as in \cite[section 3]{Th} (see Appendix), we obtain a smooth
variety $G$ which is simultaneously the blow-up of $G(\alpha^+)$ along
$G^+(\alpha)$ and the blow-up of $G(\alpha^-)$ along $G^-(\alpha)$.
Moreover the exceptional divisors of the blow-ups coincide. If we write
$$
S(n_1,d_1)=G_L(n_1,d_1,n_1-1)\times G_L(n_2,d_2,n_2-1),
$$ 
then the
exceptional divisor $Y$ is the disjoint union
$$
Y=\bigsqcup G^+(n_1,d_1)\times_{S(n_1,d_1)} G^-(n_1,d_1)
$$
for the values of $n_1$, $d_1$ which correspond to the critical value 
$\alpha$ (Appendix, (\ref{eqn:ap17})).

Now write for convenience
$$
P_{r,e}=P(G_L(r,e,r-1)).
$$

\begin{proposition} Let $e\ge\max\{1,r-g\}$. Then
\begin{equation}\label{pre}
P_{r,e}(t)=\frac{(1+t)^{2g}(1-t^{2(e+g-r+1)})(1-t^{2(e+g-r+2)})\cdots(1-
t^{2(e+g-1)})}
{(1-t^2)(1-t^4)\cdots(1-t^{2(r-1)})}.
\end{equation}\end{proposition}
\pf This follows from Proposition \ref{poincare} and (\ref{grass}).\qed
\begin{theorem}\label{poincare2}
For any
non-critical value $\alpha'$ in the interval
$(\alpha_T,\frac{d}2)$,
\begin{equation}\label{eq:pt}
 P(G(\alpha'))(t)-P(G_L)(t)=\sum\frac{t^{2C_{21}(n_1,d_1)}-
 t^{2C_{12}(n_1,d_1)}}{1-t^2}P_{n_1,d_1}(t)P_{n_2,d_2}(t),
\end{equation}
where the summation is over all solutions of {\em(\ref{eq:n1})},
{\em(\ref{eq:alphan-2})}, {\em(\ref{eq:d1})} for which $\alpha>\alpha'$. In
particular
$$
P(G_T)(t)-P(G_L)(t)=\sum\frac{t^{2C_{21}(n_1,d_1)}-t^{2C_{12}(n_1,d_1)}}{1-
t^2}P_{n_1,d_1}(t)P_{n_2,d_2}(t),
$$
the summation being over all solutions of {\em(\ref{eq:n1})},
{\em(\ref{eq:alphan-2})}, {\em(\ref{eq:d1})} for which $\alpha>\alpha_T$.
\end{theorem}

\pf We have, by Lemma \ref{disjoint}
\begin{eqnarray*}
 P(G^+(n_1,d_1))(t)&=&
 P_{n_1,d_1}(t)P_{n_2,d_2}(t)\left(1+t^2+\ldots+t^{2(C_{21}(n_1,d_1)-
1)}\right)\\
 &=&\frac{1-t^{2C_{21}(n_1,d_1)}}{1-t^2}P_{n_1,d_1}(t)P_{n_2,d_2}(t),
\end{eqnarray*}
with a similar formula for $P(G^-(n_1,d_1))(t)$. The  formula for the
Poincar\'e polynomial of a blow-up (see \cite[p.605]{GH}) now  gives
\begin{eqnarray*}
P(G)(t)&=&P(G(\alpha^+))(t)+\sum\left(t^2+\ldots +
t^{2(C_{12}(n_1,d_1)-1)}\right)P(G^+(n_1,d_1))(t)\\
&=&P(G(\alpha^-))(t)+\sum\left(t^2+\ldots +
t^{2(C_{21}(n_1,d_1)-1)}\right)P(G^-(n_1,d_1))(t),
\end{eqnarray*}
where the summation is over all solutions $(n_1,d_1)$ of
(\ref{eq:n1}), (\ref{eq:alphan-2}), (\ref{eq:d1}) for the given
$\alpha$. With a little manipulation, this gives
\be\label{eq:local}
 P(G(\alpha^-))(t)=P(G(\alpha^+))(t)+\sum
 \frac{t^{2C_{21}(n_1,d_1)}-t^{2C_{12}(n_1,d_1)}}{1-
t^2}P_{n_1,d_1}(t)P_{n_2,d_2}(t).
\end{equation}
The theorem now follows by adding the formulae (\ref{eq:local})
for all relevant $\alpha$.\qed

\begin{remark}\label{rmk7}
{\em So far, this does not depend on any coprimality assumptions, since
the flip loci lie strictly inside the moduli spaces of
$\alpha^\pm$-stable coherent systems. When $d$ is odd, however, we
can say a bit more. In this case, we know by Proposition \ref{prop:restate}(iv)
that $G_L$ is a fibration over $M(2,d)$ with fibre
$\Gr(n-2,d+2g-2)$, so we can write down an explicit formula for
$P(G_L)$. Hence (\ref{eq:pt}) gives an explicit formula for
$P(G(\alpha'))$ for $\alpha'\in(\alpha_T,\frac{d}2)$.}
\end{remark}

\begin{corollary}\label{cor:n=3}

Suppose $n=3$ and $d$ is odd. Then, for $\max\{\frac{d-3}2,0\}<\alpha'<
\frac{d}2$,
\begin{eqnarray*}
P(G(\alpha'))(t)&=&P(G_L)(t)\\&=&P(M(2,d))(t)\frac{1-t^{2(d+2g-2)}}{1-t^2}\\
&=&\frac{(1+t)^{2g}((1+t^3)^{2g}-t^{2g}(1+t)^{2g})(1-t^{2(d+2g-2)})}{(1-
t^2)^2(1-t^4)}.
\end{eqnarray*}
\end{corollary}

\pf In this case, there are no solutions to (\ref{eq:n1}) and
(\ref{eq:d1}), so $G_T=G_L$ by Theorem \ref{poincare2}. Now use Proposition
\ref{poincare} and (\ref{moduli}).\qed

\begin{corollary}\label{cor:n=4}
Suppose $n=4$ and $d$ is odd. Then
\begin{itemize}
\item[\rm(i)] if $\max\{\frac{d-2}2,0\}<\alpha'<\frac{d}2$,
\begin{eqnarray*}P(G(\alpha'))(t)&=&P(G_L)(t)\\&=&
\frac{(1+t)^{2g}((1+t^3)^{2g}-t^{2g}(1+t)^{2g})(1-t^{2(d+2g-3)})
(1-t^{2(d+2g-2)})}{(1-t^2)^2(1-t^4)^2};
\end{eqnarray*}
\item[\rm(ii)] if $\max\{\frac{d-4}2,0\}<\alpha'<\frac{d-2}2$,
\begin{eqnarray*}
P(G(\alpha'))(t)&=&
P(G_L)(t)+\frac{t^{2g}-t^{6g+2d-10}}{1-t^2}P_{1,\frac{d-1}2}(t)P_{3,\frac{d+1}2}(t)\\
&=&\frac{(1+t)^{2g}((1+t^3)^{2g}-t^{2g}(1+t)^{2g})
(1-t^{2(d+2g-3)})(1-t^{2(d+2g-2)})}{(1-t^2)^2(1-t^4)^2}\\
&&+\frac{(t^{2g}-t^{6g+2d-10})(1-t^{d-3+2g})
(1-t^{d-1+2g})(1+t)^{4g}}{(1-t^2)^2(1-t^4)}.
\end{eqnarray*}
\end{itemize}\end{corollary}

\pf For $n=4$, there is no solution to (\ref{eq:n1}), (\ref{eq:alphan-2})
and (\ref{eq:d1}) for $d=1$, but for $d\ge3$ there
is a unique solution given by
$$
n_1=1,\ \ d_1=\frac{d-1}2,\ \ \alpha=\frac{d-2}2.
$$
Moreover (\ref{C211}) and (\ref{C121}) now give
$$
C_{21}(n_1,d_1)=g,\ \ C_{12}(n_1,d_1)=3g+d-5.
$$
So, if $\alpha'>\frac{d-2}2$, $G(\alpha')=G_L$, while,
if $\alpha'\in(\frac{d-4}2,\frac{d-2}2)$, then, by Theorem \ref{poincare2},
$$
P(G(\alpha'))(t)=
P(G_L)(t)+\frac{t^{2g}-t^{6g+2d-10}}{1-t^2}P_{1,\frac{d-1}2}(t)P_{3,\frac{d+1}2}(t).
$$
The result now follows from Proposition \ref{poincare}, (\ref{moduli}), [\ref{grass}) and (\ref{pre}).\qed

\begin{remark}\label{rmk8}
{\em
The above formulae are
independent of the coefficient field for the homology groups. It
follows that the integral homology of $G(\alpha')$ is also
torsion-free. A further observation is that we can use other
cohomology theories and get similar results (for example,
algebraic cohomology). With some further work, it may also be
possible to compute Chow groups.}
\end{remark}

\section*{Appendix. Geometry of flips}

\renewcommand{\thesection}{A}
\setcounter{theorem}{0}

Our object in this appendix is to establish the geometric
description of a flip in the best-behaved case. This is analogous
to Thaddeus' description \cite{Th} for the case $n=2$, $k=1$. The
general approach is the same as that of Thaddeus and is 
similar to that of He \cite{He} (see also \cite{Sep} and, for a related
problem, \cite{GGM}). However, He restricts to rank $2$
(and also to coherent sheaves on $\PP^2$) at a crucial point,
although earlier in the paper he discusses coherent systems in
great generality.

As in the main part of the paper, we work with coherent systems
over a smooth projective irreducible algebraic curve $C$, defined
over the complex numbers. For a fixed type $(n,d,k)$, let
$G(\alpha)=G(\alpha; n,d,k)$ and let $\alpha_c$ be a critical
value of $\alpha$. The flip locus $G^+:=G^+(\alpha_c)\subset G(\alpha_c^+)$ is given by
non-trivial extensions of the form
 \begin{equation}\label{eqn:ap1}
  0 \lra (E_1,V_1) \lra (E,V) \lra (E_2,V_2)\lra 0,
 \end{equation}
where $(E_j,V_j)$ is of type $(n_j,d_j,k_j)$ and is
$\alpha_c$-semistable and $\alpha_c^+$-stable for $j=1,2$.
Moreover
 \begin{equation}\label{eqn:ap2}
  \mu_{\alpha_c}(E_1,V_1)=\mu_{\alpha_c}(E_2,V_2), \quad
  \frac{k_1}{n_1}<\frac{k_2}{n_2}.
 \end{equation}
In particular, we have
 $$
 \alpha_c=\frac{n_2d_1-n_1d_2}{n_1k_2-n_2k_1}>0,
 $$
hence $\frac{d_1}{n_1}>\frac{d_2}{n_2}$. The flip locus
$G^-:=G^-(\alpha_c)\subset G(\alpha_c^-)$ is given similarly by extensions
 \begin{equation}\label{eqn:ap3}
  0 \lra (E_2,V_2) \lra (E',V') \lra (E_1,V_1)\lra 0,
 \end{equation}
where $(E_j,V_j)$ is now $\alpha_c^-$-stable for $j=1,2$, the
other conditions being as above. Now $G(\alpha_c^-)$ is obtained
from $G(\alpha_c^+)$ by deleting $G^+$ and inserting $G^-$; in
particular
 $$
 G(\alpha_c^-) \setminus G^- = G(\alpha_c^+) \setminus G^+ \, .
 $$

It is possible that there is more than one choice of the values
$n_1, d_1,k_1$ for a given critical value $\alpha_c$. In this case
we write $G^+(n_1,d_1)$, $G^-(n_1,d_1)$ for the subsets of the
flip loci corresponding to particular values of $n_1,d_1$; note
that, once $n_1,d_1$ are fixed, so is $k_1$.

\begin{assumption} \label{ass:ap1}
We assume, for all choices of $n_1,d_1,k_1$ corresponding to the
critical value $\alpha_c$,
 \begin{enumerate}
 \item[(a)] $\GCD (n_1,d_1,k_1) = \GCD (n_2,d_2,k_2) =1$\ ;
 \item[(b)] $\alpha_c$ is not a critical value for $(n_1,d_1,k_1)$,
 $(n_2,d_2,k_2)$\ ;
 \item[(c)] $G_1:=G(\alpha_c;n_1,d_1,k_1)$ and $G_2:=G(\alpha_c;n_2,d_2,k_2)$
 are smooth of the expected dimensions\ ;
 \item[(d)] $\Ext^2((E_1,V_1),(E_2,V_2))=
 \Ext^2((E_2,V_2),(E_1,V_1))=0$.
 \end{enumerate}
\end{assumption}

\begin{remark}\begin{em}
Note that, by Theorem \ref{th:inj} and the proof of 
Lemma \ref{lemma:h0=h2=0}, (c) and (d) 
always hold if $\alpha_c$ is in the injective range for both $(n_1,d_1,k_1)$
and $(n_2,d_2,k_2)$. 
\end{em}\end{remark}

\begin{lemma} \label{lem:ap2}
 Given Assumptions \ref{ass:ap1}, we have
 \begin{enumerate}
 \item[(i)] $(E_j,V_j) \in G_j$ for $j=1,2$\ ;
 \item[(ii)] $G_1$, $G_2$ are smooth projective varieties\ ;
 \item[(iii)] $\Hom ((E_1,V_1),(E_2,V_2))=
 \Hom((E_2,V_2),(E_1,V_1))=0$.
 \end{enumerate}
\end{lemma}

\pf (i) follows at once from (a) and (b). (ii) follows from (a) and (c).
Finally, for (iii), it follows from (\ref{eqn:ap2}) that
$(E_1,V_1)$, $(E_2,V_2)$ are non-isomorphic with the same
$\alpha_c$-slope; since both are $\alpha_c$-stable by (i), this
implies (iii). \qed

\begin{corollary} \label{cor:ap3}
 The sequences {\rm (\ref{eqn:ap1})}, {\rm (\ref{eqn:ap3})} are the unique
 Jordan-H\"older filtrations of $(E,V)$, $(E',V')$ as
 $\alpha_c$-semistable coherent systems.
\end{corollary}

\pf
 This follows from (i) and the non-triviality of (\ref{eqn:ap1}),
 (\ref{eqn:ap3}).
\qed

Our next object is to show that $G^+$, $G^-$ are smooth
subvarieties of $G(\alpha_c^+)$, $G(\alpha_c^-)$ respectively, and
to identify their normal bundles. The proofs are identical in the
two cases, so we shall work with the $+$ case only. Since we have
not assumed that $G(\alpha_c^+)$ is smooth, we need first a lemma.

\begin{lemma} \label{lem:ap4}
 $G(\alpha_c^+)$ is smooth of the expected dimension at every point
 of $G^+$.
\end{lemma}

\pf
 We must show that, under Assumptions \ref{ass:ap1}, $G(\alpha_c^+)$
 is smooth of the expected dimension at $(E,V)$ for any non-trivial
 extension (\ref{eqn:ap1}). For this, we must show that the Petri
 map is injective, or equivalently that
   $$
   H^0(E^*\otimes N\otimes K)=0 ,
   $$
 where $N$ is the kernel of the evaluation map $V\otimes \cO \lra
 E$. Let $N_l$ be the kernel of the evaluation map
 $V_l\otimes \cO \lra E_l$, for $l=1,2$. Then
 $H^0(E_m^* \otimes N_l\otimes K)=0$ for $l=1,2$, $m=1,2$,
 by Assumption \ref{ass:ap1} (c) and (d), using \cite[Proposition 3.2]{BGMN}.
 The result now follows by diagram chasing. \qed

\begin{definition}\label{def:ap1}
{\em A {\em family of coherent systems of type $(n,d,k)$ on $C$
parametrised by $S$} is a pair $(\cE,\cV)$, where $\cE$ is a vector
bundle of rank $n$ over $S\times C$ such that
$\cE_s=\cE|_{\{s\}\times C}$ has degree $d$ for all $s\in S$, and
$\cV$ is a locally free subsheaf of $p_{S*} \cE $ (where $p_S: S
\times C \lra S$ stands for the projection) of rank $k$ such that
the fibres $\cV_s$ map injectively to $H^0(\cE_s)$ for all $s\in
S$.}
\end{definition} 

\begin{remark}
{\em This definition is more restrictive than that of \cite{He}, but is 
sufficient for our purposes.}
\end{remark}

We need two facts about families of coherent systems
for which we have been unable to locate proofs. We state these as propositions.

\begin{proposition}\label{prop:ap1} Suppose $\GCD(n,d,k)=1$. Then there exists a 
universal family of coherent systems over $G(\alpha;n,d,k)\times C$.
\end{proposition}

\pf 
[The proof is on standard lines and is modelled on those of \cite{Ty2} for 
$M(n,d)$ (see also \cite[Theorem 5.12]{N2} or \cite[Premi\`ere Partie, 
Th\'eor\`eme 18]{Ses}) when $\GCD(n,d)=1$ and \cite[Theorems 2.8 and 3.3]{RV} 
for $G(\alpha;n,d,k)$ when $\alpha$ is small and $\GCD(n,k)=1$.]

We recall the method of construction of $G(\alpha;n,d,k)$ \cite{LeP,RV,KN}.
We have a family $(\cE,\cV)$ of $\alpha$-stable coherent systems of 
type $(n,d,k)$ parametrised by a variety $R^s$ together with an action of
$PGL(N)$ on $R^s$ which lifts to an action of $GL(N)$ on $(\cE,\cV)$.
The family $(\cE,\cV)$ satisfies the local universal property for
$\alpha$-stable coherent systems of type $(n,d,k)$ and the moduli space
$G(\alpha;n,d,k)$ is the geometric quotient of $R^s$ by $PGL(N)$ (indeed
$R^s$ is a principal $PGL(N)$-fibration over $G(\alpha;n,d,k)$). Moreover
the action of an element $\lambda$ of the centre $\CC^*$ of $GL(N)$ on 
$(\cE,\cV)$ is multiplication by $\lambda$. Suppose now that we can construct
a line bundle $\cL$ on $R^s$ such that the action of $PGL(N)$ on $R^s$ lifts
to an action of $GL(N)$ on $\cL$ with the same property. We then consider
the coherent system
$$
(\cE\otimes p^*_{R^s}\cL^*,\cV\otimes\cL^*)
$$
over $R^s\times C$. The action of $PGL(N)$ on $R^s\times C$ now lifts to an action of
$PGL(N)$ on this coherent system. It follows from the theory of descent
(see \cite[Theorem 1]{Gr} or Kempf's descent lemma \cite[Theorem 2.3]{DN}) 
that this coherent system is the pull-back of a coherent system over
$G(\alpha;n,d,k)\times C$ which satisfies the required universal property.  

It remains to construct $\cL$. For this, we consider the bundle $\cE$ over 
$R^s\times C$ and let $\cE_t$ denote the bundle over $C$ obtained by
restricting $\cE$ to $\{t\}\times C$. There exists a line bundle $L$ on
$C$ such that $H^1(\cE_t\otimes L)=0$ for all $t\in R^s$. It follows by
Riemann-Roch that
$$
p:=h^0(\cE_t\otimes L)=d+n(m+1-g)
$$
for all $t\in R^s$, where $m=\deg L$. Hence
$$
\cF:=p_{R^s*}(\cE\otimes p^*_CL)
$$
is locally free of rank $p$.

Now choose a point $x_0\in C$; then, by the same argument,
$$
\cF(x_0):=p_{R^s*}(\cE\otimes p^*_CL(x_0))
$$
is locally free of rank $p+n$. Now
$$
\GCD(p+n,p,k)=\GCD(n,p,k)=\GCD(n,d,k)=1,
$$
so there exist integers $a$, $b$, $c$ such that
$$
a(p+n)+bp+ck=1.
$$
We can now define
$$
\cL:=(\det\cF(x_0))^a\otimes(\det\cF)^b\otimes(\det\cV)^c.
$$
The element $\lambda\in\CC^*$ now acts on $\cL$ by
$$
\lambda^{a(p+n)+bp+ck}=\lambda
$$
as required.
\qed

\begin{proposition}\label{prop:ap2} Let $(\cE_1,\cV_1)$, $(\cE_2,\cV_2)$
be two families of coherent systems parametrised by $S$ and let
 $$
   \cExt^q_{p_S} ((\cE_2,\cV_2),(\cE_1,\cV_1))
  $$
be defined as in \cite[1.2]{He}. Then there exists a spectral sequence 
with $E_2$-term
$$
E^{pq}_2=H^p(\cExt^q_{p_S} ((\cE_2,\cV_2),(\cE_1,\cV_1))),
$$
which abuts to $\Ext^*((\cE_2,\cV_2),(\cE_1,\cV_1))$.
\end{proposition}

\pf The construction of \cite{He} depends on embedding the category of 
families of coherent systems over $C$ parametrised by $S$ into an abelian 
category $\cC$ with enough injectives; in He's notation, the objects of this 
larger
category $\cC$ are called {\em algebraic systems on $S\times X$ relative to
$S$}. Both $\cExt_{p_S}$ and $\Ext$ can now be defined using an injective
resolution of $(\cE_1,\cV_1)$ in $\cC$. To prove the existence of the
spectral sequence, it is sufficient to show that, if $\cI$ is an injective 
in $\cC$, then $\cHom_{p_S}((\cE_2,\cV_2),\cI)$ is an acyclic sheaf;
this follows from He's description of the injectives \cite[Th\'eor\`eme 1.3]{He}.
\qed

\begin{lemma} \label{lem:ap5}
There exists a vector bundle $W^+$ over $G_1\times G_2$ and a
morphism
  $$
  f_+: \PP W^+ \lra G(\alpha_c^+),
  $$
which maps $\PP W^+$ bijectively to $G^+(n_1,d_1)$.
\end{lemma}

\pf By Assumption \ref{ass:ap1} (a) and Proposition \ref{prop:ap1},
there exist
universal families of coherent systems $(\cE_1,\cV_1)$ over
$G_1\times C$ and $(\cE_2,\cV_2)$ over $G_2\times C$.
By Assumption \ref{ass:ap1} (d) and Lemma \ref{lem:ap2} (iii),
  $$
  \dim \Ext^1((E_2,V_2),(E_1,V_1))
  $$
is independent of the choice of $(E_1,V_1) \in G_1$, $(E_2,V_2) \in
G_2$. It follows from \cite[Corollaire 1.20]{He} that there is a vector bundle
$W^+$ over $G_1\times G_2$ whose fibre over 
$$((E_1,V_1),(E_2,V_2)) \in G_1\times G_2
$$
is $\Ext^1((E_2,V_2),(E_1,V_1))$; in fact, in the notation of
\cite{He},
  $$
  W^+= \cExt^1_\pi ((p_2\times \id)^*(\cE_2,\cV_2),
  (p_1\times \id)^*(\cE_1,\cV_1)),
  $$
where $\pi:G_1\times G_2\times C \lra G_1\times G_2$,
$p_1:G_1\times G_2 \lra G_1$ and $p_2:G_1\times G_2 \lra G_2$ are
the natural projections. Now $\PP W^+$ classifies the non-trivial
extensions (\ref{eqn:ap1}) up to scalar multiples. We can
therefore define $f_+$ set-theoretically as the natural map
sending (\ref{eqn:ap1}) to $(E,V)\in G(\alpha_c^+)$. The fact that
$f_+$ maps $\PP W^+$ bijectively to $G^+(n_1,d_1)$ follows from
Corollary \ref{cor:ap3}.

In order to prove that $f_+$ is a morphism, we need to construct a
universal extension (\ref{eqn:ap1}) over $\PP W^+\times C$. This
is done in exactly the same way as for extensions of vector
bundles. Let $\sigma: \PP W^+\times C \lra \PP W^+$ and $p: \PP
W^+\lra G_1\times G_2$ be the natural projections. We write, for
$m=1,2$,
  $$
  (\cE_m, \cV_m)^+ =(p\times \id)^* (p_m\times \id)^* (\cE_m,\cV_m).
  $$
We construct the universal extension as an extension
 \begin{equation}\label{eqn:ap4}
  0 \lra (\cE_1, \cV_1)^+ \otimes \sigma^* \cO_{\PP W^+}(1) \lra
  (\cE, \cV)^+ \lra (\cE_2, \cV_2)^+ \lra 0
 \end{equation}
on $\PP W^+\times X$. Extensions of the form (\ref{eqn:ap4}) are
classified by
 $$
 \Ext^1((\cE_2, \cV_2)^+,(\cE_1, \cV_1)^+ \otimes \sigma^* \cO_{\PP
 W^+}(1)).
 $$
By Proposition \ref{prop:ap2} we have a spectral sequence whose
$E_2$-term is given by
  $$
  E_2^{pq}=H^p \big(\cExt^q_\sigma ((\cE_2, \cV_2)^+,
  (\cE_1, \cV_1)^+ \otimes \sigma^* \cO_{\PP W^+}(1))\big)
  $$
which abuts to $\Ext^*((\cE_2, \cV_2)^+,(\cE_1, \cV_1)^+ \otimes
\sigma^* \cO_{\PP W^+}(1))$. In view of Assumption \ref{ass:ap1}
(d) and Lemma \ref{lem:ap2} (iii), we have $E_2^{pq}=0$ except for
$q=1$. Hence
 $$
  \begin{aligned}
  \Ext^1 &((\cE_2, \cV_2)^+,(\cE_1, \cV_1)^+ \otimes \sigma^* \cO_{\PP
  W^+}(1)) \\ &\cong H^0\big(\cExt^1_\sigma ((\cE_2, \cV_2)^+,
  (\cE_1, \cV_1)^+ \otimes \sigma^* \cO_{\PP W^+}(1))\big).
  \end{aligned}
  $$
On the other hand, by base-change \cite[Th\'eor\`eme 1.16]{He},
  $$
  \begin{aligned}
  \cExt^1_\sigma &((\cE_2, \cV_2)^+,
  (\cE_1, \cV_1)^+ \otimes \sigma^* \cO_{\PP W^+}(1))
  \\
  &\cong p^* \cExt_\pi^1 ((p_2\times \id)^* (\cE_2,\cV_2),
  (p_1\times \id)^* (\cE_1,\cV_1)) \otimes \cO_{\PP W^+}(1) \\ &=
  p^* W^+ \otimes \cO_{\PP W^+}(1).
  \end{aligned}
  $$
Now $H^0(p^* W^+ \otimes \cO_{\PP W^+}(1))= \End W^+$. The
universal extension is then the extension (\ref{eqn:ap4})
corresponding to the identity endomorphism of $W^+$. It is clear
that the restriction of (\ref{eqn:ap4}) to $\{y\}\times X$ is
precisely the extension (\ref{eqn:ap1}) corresponding to $y\in \PP
W^+$. So the morphism $\PP W^+ \lra G(\alpha_c^+)$, given by
(\ref{eqn:ap4}) and the universal property of $G(\alpha_c^+)$,
coincides with $f_+$. This completes the proof of the lemma. \qed

It is an immediate consequence of this lemma and Lemma
\ref{lem:ap4} that $G^+(n_1,d_1)$ is a projective subvariety of
the smooth part of $G(\alpha_c^+)$. Moreover it follows from
Corollary \ref{cor:ap3} that the $G^+(n_1,d_1)$ for different
values of $n_1$ and $d_1$  are disjoint. The
computation of the normal bundle can therefore be carried out
independently for each choice of $n_1,d_1$. Let $W^-$ be the
bundle over $G_1\times G_2$ constructed in an analogous way to
$W^+$ after interchanging the subscripts $1,2$.

\begin{proposition}\label{prop:ap6}
 The morphism $f_+$ is a smooth embedding with normal bundle $p^*W^-
 \otimes \cO_{\PP W^+}(-1)$.
\end{proposition}

\pf The proof is exactly analogous to \cite[(3.9)]{Th}. In our
notation, it proceeds as follows. For simplicity, we begin by
looking at the infinitesimal deformations at a point $\xi$ of $\PP
W^+$ represented by an extension (\ref{eqn:ap1}). We have a short
exact sequence of complexes
  $$
  \begin{array}{ccccccccc}
  0 & \lra & A & \lra & \Hom(E,E) &\lra & \Hom (E_1,E_2) &\lra &0\\
    && \downarrow &&\downarrow && \downarrow \\
  0 & \lra & B & \lra & \Hom(V,E/V) &\lra & \Hom (V_1,E_2/V_2) &\lra
  &0.
  \end{array}
  $$
Here $V, V_1,V_2$ are to be interpreted as sheaves of locally
constant sections, the right-hand square is the obvious
homomorphism of complexes and $A\to B$ is just the kernel of this
homomorphism. The middle complex parametrises infinitesimal
deformations of $(E,V)$, while the right-hand one is the complex
giving rise to the fibre of $W^-$ at the point $((E_1,V_1),
(E_2,V_2))$ of $G_1\times G_2$ (see the proof of Lemma
\ref{lem:ap5} and \cite[Corollaire 1.6]{He}). Taking hypercohomology, we
therefore obtain an exact sequence
  \begin{equation}\label{eqn:ap14}
  0\lra \HH^1(A\to B) \lra TG(\alpha_c^+)_{f_+(\xi)} \lra
  W^-_{p(\xi)} \lra 0.
  \end{equation}
Now we have another natural short exact sequence of complexes
  $$
  \begin{array}{ccccccccc}
  0 & \lra & \Hom(E_2,E_1) & \lra & A &\lra & \Hom (E_1,E_1) \oplus \Hom(E_2,E_2) &\lra &0\\
    && \downarrow &&\downarrow && \downarrow \\
  0 & \lra & \Hom(V_2,E_1/V_1) & \lra & B &\lra & \Hom (V_1,E_1/V_1) \oplus \Hom(V_2,E_2/V_2) &\lra
  &0.
  \end{array}
  $$
In this case the hypercohomology gives
  $$
  0 \lra (T_{fibre} \PP W^+)_{\xi} \lra \HH^1(A\to B) \lra T(G_1\times
  G_2)_{p(\xi)} \lra 0.
  $$
This sequence identifies $\HH^1(A\to B)$ with $(T\PP W^+)_\xi$.
The sequence (\ref{eqn:ap14}) now becomes
  \begin{equation}\label{eqn:ap15}
  0 \lra (T\PP W^+)_\xi \stackrel{df_+}{\lra} T
  G(\alpha_c^+)_{f_+(\xi)}
  \lra W^-_{p(\xi)} \lra 0.
  \end{equation}
Since we know that $G(\alpha_c^+)$ is smooth at $f_+(\xi)$, this
shows that $f_+$ is smooth at $\xi$ and identifies the normal
space with $W^-_{p(\xi)}$.

It remains to globalise this construction. For this we need to
replace $(E,V)$, $(E_1,V_1)$, $(E_2,V_2)$ by $(\cE, \cV)^+$,
$(\cE_1, \cV_1)^+\otimes  \sigma^*\cO_{\PP W^+}(1)$, $(\cE_2,
\cV_2)^+$ in accordance with (\ref{eqn:ap4}). The sequence
(\ref{eqn:ap15}) now becomes
  $$
 0 \lra T\PP W^+ \lra f_+^* T G(\alpha_c^+)
  \lra p^*W^- \otimes \cO_{\PP W^+} (-1) \lra 0.
 $$
Since we already know from Lemma \ref{lem:ap5} that $f_+$ is
injective, this completes the proof of the proposition. \qed

If we now blow up $G(\alpha_c^+)$ along $G^+(n_1,d_1)$, it follows
from Proposition \ref{prop:ap6} that the exceptional divisor is
isomorphic to $\PP W^+ \times_{G_1\times G_2} \PP W^-$. This works
for each allowable choice of $n_1,d_1$. To avoid confusion, we
label these choices by $1,\ldots, r$ and denote the corresponding
$W^+,W^-$ by $W^+_j, W^-_j$ for $1\leq j \leq r$. Now performing
all the blow-ups simultaneously, we obtain a variety
$\widetilde{G(\alpha_c^+)}$ with exceptional divisors $Y_1,
\ldots, Y_r$, all contained in the smooth part of
$\widetilde{G(\alpha_c^+)}$. In exactly the same way, we blow-up
$G(\alpha_c^-)$ along the various $G^-(n_1,d_1)$ (labelled as
before) to obtain $\widetilde{G(\alpha_c^-)}$ with exceptional
divisors $Y_1',\ldots, Y_r'$. Note that there exist natural
isomorphisms
\begin{equation}\label{eqn:ap17}
  Y_j \cong \PP W_j^+ \times_{G_1\times G_2} \PP W_j^- \cong Y_j' .
\end{equation}
We use these isomorphisms to identify $Y_j$ and $Y_j'$.

The final step in the construction is to show that
$\widetilde{G(\alpha_c^+)}$ is naturally isomorphic to
$\widetilde{G(\alpha_c^-)}$. It is easy to construct a natural
bijection between these varieties. In fact, if we write $Y=Y_1\cup
\cdots \cup Y_r$ and $Y'=Y_1'\cup \cdots \cup Y_r'$ then
  $$
  \widetilde{G(\alpha_c^+)} \setminus Y =
  \widetilde{G(\alpha_c^-)}\setminus Y',
  $$
each variety consisting precisely of the $\alpha_c$-stable
coherent systems. On the other hand, as observed above, $Y_j$ and
$Y_j'$ can be identified. It remains to show that there exist
morphisms $\widetilde{G(\alpha_c^+)} \lra
\widetilde{G(\alpha_c^-)}$ and $\widetilde{G(\alpha_c^-)} \lra
\widetilde{G(\alpha_c^+)}$ such that the following diagram
commutes for each $j$\ :
  \begin{equation}\label{eqn:ap5}
  \begin{array}{ccccc}
  \widetilde{G(\alpha_c^+)} \setminus Y & \subset &
  \widetilde{G(\alpha_c^+)} & \hookleftarrow & Y_j \\
   \parallel & & \downarrow\uparrow && \parallel \\
  \widetilde{G(\alpha_c^-)} \setminus Y' & \subset &
  \widetilde{G(\alpha_c^-)} &\hookleftarrow & Y_j'.
  \end{array}
  \end{equation}

For this purpose, we prove

\begin{proposition}\label{prop:ap7}
 There exists a morphism $\widetilde{G(\alpha_c^+)} \lra
 G(\alpha_c^-)$ making the following diagram commute:
  \begin{equation}\label{eqn:ap6}
  \begin{array}{ccccc}
  \widetilde{G(\alpha_c^+)} \setminus Y & \subset &
  \widetilde{G(\alpha_c^+)} & \hookleftarrow & Y_j \\
   \downarrow & & \downarrow && \ \downarrow q\\
  G(\alpha_c^-) \setminus (\PP W_1^- \cup \cdots \cup \PP W_r^-) & \subset &
  G(\alpha_c^-) & \hookleftarrow & \PP W_j^-,
  \end{array}
  \end{equation}
where $q$ is the natural projection.
\end{proposition}

\pf
 Suppose first that $\GCD(n,d,k)=1$, so that there exists a
 universal coherent system on $G(\alpha_c^+)\times C$. We write
 $(\cE,\cV)$ for the pull-back of such a coherent system to
 $\widetilde{G(\alpha_c^+)}\times C$. We want to compare
 $(\cE,\cV)|_{Y_j\times C}$ with the pull-backs of the extension
 (\ref{eqn:ap4}) and the equivalent extension for $\PP
 W_j^-\times C$ to $Y_j\times C$.

 Let $r_j^+:Y_j \lra \PP W_j^+$ and $r_j^-:Y_j \lra \PP W_j^-$
 denote the projections.
 We write (\ref{eqn:ap4}) for $\PP W_j^+\times C$ as
 \begin{equation}\label{eqn:ap7}
  0 \lra (\cE_1, \cV_1)^+_j \otimes (\sigma^+_j)^* \cO_{\PP W^+_j}(1) \lra
  (\cE, \cV)^+_j \lra (\cE_2, \cV_2)^+_j \lra 0,
 \end{equation}
 where  $\sigma_j^+:\PP W_j^+\times C \lra  \PP W_j^+$,
 and the corresponding extension on $\PP W_j^-\times C$ as
 \begin{equation}\label{eqn:ap8}
  0 \lra (\cE_2, \cV_2)^-_j \otimes (\sigma^-_j)^* \cO_{\PP W^-_j}(1) \lra
  (\cE, \cV)^-_j \lra (\cE_1, \cV_1)^-_j \lra 0,
 \end{equation}
 where $\sigma_j^-:\PP W_j^-\times C \lra \PP W_j^-$.
 Since $(\cE_1,\cV_1)^\pm_j$ are both pulled back from the same
 coherent system on $G_1\times G_2\times C$, we have
 \begin{equation}\label{eqn:ap9}
  (r_j^+\times \id)^*(\cE_1,\cV_1)_j^+ =  (r_j^-\times
  \id)^*(\cE_1,\cV_1)_j^-,
 \end{equation}
 with a similar statement for $(\cE_2,\cV_2)_j^\pm$. For every $y\in
 Y_j$, the restrictions of $(\cE,\cV)$ and $(r_j^+\times
 \id)^*(\cE,\cV)_j^+$ to $\{y\}\times C$ are isomorphic and
 the
 coherent systems are all $\alpha_c^+$-stable; it follows that
 $$
  (\cE,\cV)|_{Y_j\times C} \cong (r_j^+\times \id)^* (\cE,
  \cV)^+_j \otimes L_j
 $$
 for some line bundle $L_j$ pulled back from $Y_j$. We define
 $$
   \widetilde{(\cE_1,\cV_1)}_j := (r_j^+\times \id)^*
   (\cE_1,\cV_1)_j^+ \otimes L_j,
 $$
 with a similar definition for $\widetilde{(\cE_2,\cV_2)}_j$. We
 write also
 $$
 \cO_{Y_j\times C}(a,b) = (r_j^+\times
 \id)^*(\sigma_j^+)^*\cO_{\PP W_j^+}(a)\otimes
 (r_j^-\times
 \id)^*(\sigma_j^-)^*\cO_{\PP W_j^-}(b).
 $$
 Taking account of (\ref{eqn:ap9}), we can then
 tensor
 the pull-backs of (\ref{eqn:ap7}) and (\ref{eqn:ap8}) to
 $Y_j\times C$ by $L_j$ to get
 \begin{equation}\label{eqn:ap10}
 0\lra \widetilde{(\cE_1,\cV_1)}_j \otimes \cO_{Y_j\times C}(1,0)
 \lra  (\cE,\cV)|_{Y_j\times C} \lra \widetilde{(\cE_2,\cV_2)}_j
 \lra 0
 \end{equation}
 and
 \begin{equation}\label{eqn:ap11}
 0\lra \widetilde{(\cE_2,\cV_2)}_j \otimes \cO_{Y_j\times C}(0,1)
 \lra  (r_j^- \times \id)^*(\cE,\cV)_j^- \otimes L_j
 \lra \widetilde{(\cE_1,\cV_1)}_j \lra 0.
 \end{equation}

 Now let
  $$
  \widetilde{(\cE_1,\cV_1)} =\bigsqcup_{j=1}^r
  \widetilde{(\cE_1,\cV_1)}_j ,\qquad \widetilde{(\cE_2,\cV_2)}
  =\bigsqcup_{j=1}^r \widetilde{(\cE_2,\cV_2)}_j,
  $$
  and define $(\hat{\cE},\hat{\cV})$ on
  $\widetilde{G(\alpha_c^+)}\times C$ by the exact sequence
  \begin{equation}\label{eqn:ap12}
   0\lra (\hat{\cE},\hat{\cV}) \lra (\cE,\cV)
   \stackrel{\psi}{\lra}\widetilde{(\cE_2,\cV_2)} \lra 0,
  \end{equation}
  where $\psi$ is given by the composition
  $$
  (\cE,\cV) \lra (\cE,\cV)|_{Y_j\times C} \lra
  \widetilde{(\cE_2,\cV_2)}_j
  $$
  in the neighbourhood of $Y_j\times C$. Since the $Y_j$ are
  disjoint Cartier divisors in $\widetilde{G(\alpha_c^+)}$,
  $\hat{\cE}$ is locally free and $(\hat{\cE},\hat{\cV})$ is a
  family of coherent systems parametrised by
  $\widetilde{G(\alpha_c^+)}$. The restriction of
  (\ref{eqn:ap12}) to $Y_j\times C$ gives a $4$-term exact
  sequence which can be split into two short exact sequences. The
  right-hand one coincides by construction with (\ref{eqn:ap10}),
  while the left-hand one takes the form
  \begin{equation}\label{eqn:ap13}
   0\lra \ker \lra (\hat{\cE},\hat{\cV})|_{Y_j\times C} \lra
   \widetilde{(\cE_1,\cV_1)}_j \otimes \cO_{Y_j\times C}(1,0) \lra
   0.
  \end{equation}
Here
\begin{eqnarray*}
\ker&\cong&\widetilde{(\cE_2,\cV_2)}_j\otimes\cTor_1(\cO_{Y_j\times C},
\cO_{Y_j\times C})\\
&\cong&\widetilde{(\cE_2,\cV_2)}_j\otimes N^*,
\end{eqnarray*}
where $N$ is the pull-back to $Y_j\times C$ of the normal bundle
of $Y_j$ in $\widetilde{G(\alpha^+_c)}$. Now, by Proposition 
\ref{prop:ap6} and a standard property of blow-ups, we have
$$
N\cong\cO_{Y_j\times C}(-1,-1).
$$
So (\ref{eqn:ap13}) becomes
\begin{equation}\label{eqn:ap16}
0\lra \widetilde{(\cE_2,\cV_2)}_j\otimes\cO_{Y_j\times C}(1,1) \lra 
(\hat{\cE},\hat{\cV})|_{Y_j\times C} \lra
   \widetilde{(\cE_1,\cV_1)}_j \otimes \cO_{Y_j\times C}(1,0) \lra
   0.
\end{equation}

In the neighbourhood of $Y_j\times C$, we can interpret this in terms 
of a commutative diagram
 $$
 \begin{array}{ccccccccc}
  && 0 && 0 \\
  && \downarrow && \downarrow \\
 && (\cE,\cV) \otimes \cO (-Y_j\times C) &=&
   (\cE,\cV) \otimes \cO (-Y_j\times C) \\
   && \downarrow && \downarrow \\
 0 & \rightarrow & (\hat{\cE},\hat{\cV}) & \rightarrow &
   (\cE,\cV) & \rightarrow &  \widetilde{(\cE_2,\cV_2)}_j & \rightarrow & 0\\
  && \downarrow && \downarrow&& || \\
 0 & \rightarrow &
   \widetilde{(\cE_1,\cV_1)}_j \otimes \cO_{Y_j\times C}(1,0) & 
   \rightarrow & (\cE,\cV)|_{Y_j\times C} & 
   \rightarrow &  \widetilde{(\cE_2,\cV_2)}_j &\rightarrow & 0\\
  && \downarrow && \downarrow \\
  && 0 && 0 
 \end{array}
 $$
Here the middle row is (\ref{eqn:ap12}), the bottom row is
(\ref{eqn:ap10}) and the quotient map in the 
left-hand column factorises as
$$
(\hat{\cE},\hat{\cV})\lra(\hat{\cE},\hat{\cV})|_{Y_j\times C} \lra
   \widetilde{(\cE_1,\cV_1)}_j \otimes \cO_{Y_j\times C}(1,0),
$$
where the two maps here are the restriction to $Y_j\times C$ and 
the quotient map of (\ref{eqn:ap16}). We are now in the situation of
the main diagram of \cite[p.583]{He}.

Using the identification (\ref{eqn:ap17}), a point of $Y_j$ can be
represented as $(t^+,t^-)$, where $t^+\in\PP W_j^+$, $t^-\in\PP W_j^-$
represent the classes of non-trivial extensions (\ref{eqn:ap1})
and (\ref{eqn:ap3}) respectively. Moreover the natural projections
to $G(\alpha_c^+)$ and $G(\alpha_c^-)$ are given by
$$
(t^+,t^-)\mapsto f_+(t^+)=(E,V),\quad(t^+,t^-)\mapsto f_-(t^-)=(E',V'),
$$
where $f_+$ is defined in Lemma \ref{lem:ap5} and $f_-$ in a similar way.
But, by Proposition \ref{prop:ap6}, $t^-$ represents a normal
direction to $G_+(n_1,d_1)$ at the point $f_+(t^+)$. He's argument
\cite[pp.583, 584]{He} now shows that $t^-$ is precisely the class
of the restriction of (\ref{eqn:ap16}) to $\{(t^+,t^-)\}\times C$.
[He doesn't quite claim this directly, but note the sentence beginning 
``Mais cette question...'' in the middle of p.584.] So, considered as 
families of extensions over $C$ parametrised by $Y_j$, the sequences
(\ref{eqn:ap16}) and (\ref{eqn:ap11}) coincide. By the universality 
of families of extensions, it follows that (\ref{eqn:ap16}) can be
obtained from (\ref{eqn:ap11}) by tensoring by $\cO_{Y_j\times C}(1,0)$. 
Hence $(\hat{\cE},\hat{\cV})$ is
a family of $\alpha_c^-$-stable coherent systems, confirming the
existence of the required morphism. The commutativity of
(\ref{eqn:ap6}) is obvious.

 If $\GCD(n,d,k)\neq 1$, we no longer have a universal family
 defined on $G(\alpha_c^+)\times C$. However
 $G(\alpha_c^+)$ is constructed as a geometric quotient
of a variety $R^s$ by an action of $PGL(N)$ such that there exists
a locally universal family on $R^s\times C$. By pulling everything back 
to $R^s$, we see that the argument above determines a blow-up 
$\widetilde{R^s}$ and a morphism
$\widetilde{R^s}\lra G(\alpha_c^-)$, which (as a map) factors through
$\widetilde{G(\alpha_c^+)}$. Since $\widetilde{R^s}\lra\widetilde{G(\alpha_c^+)}$
is again a geometric quotient, it follows that the map 
$\widetilde{G(\alpha_c^+)}\lra G(\alpha_c^-)$ is a morphism as
required. This completes the proof.
\qed

By Proposition \ref{prop:ap7}, we have a morphism
  $$
  \widetilde{G(\alpha_c^+)} \lra
  G(\alpha_c^+) \times G(\alpha_c^-),
  $$
which is injective and is easily seen to be a smooth embedding.
Moreover the image of $\widetilde{G(\alpha_c^+)}$ is precisely the
closure of the graph of the identification map
  $$
  G(\alpha_c^+) \setminus (\PP W_1^+\cup \cdots \PP W_r^+) =
  G(\alpha_c^-) \setminus (\PP W_1^-\cup \cdots \PP W_r^-)
  $$
  in $G(\alpha^+_c)\times G(\alpha^-_c)$.
By the same argument with $+$ and $-$ interchanged,
$\widetilde{G(\alpha_c^-)}$ can also be identified with the
closure of the graph of the same identification. This completes
the construction of the diagram (\ref{eqn:ap5}).

\end{document}